\documentclass[12pt]{amsart}
\usepackage{amsfonts, amsmath, amsthm}
\usepackage{graphicx}

\newtheorem{pro}{Proposition}[section]
\newtheorem{thm}[pro]{Theorem}
\newtheorem{lem}[pro]{Lemma}

\newtheorem{cnj}[pro]{Conjecture}
\newtheorem{rmkk}[pro]{Remark}
\newtheorem{rmkks}[pro]{Remarks}

\newtheorem{ass}{Assertion}
\newtheorem{clm}{Claim}

\newtheorem{cor}[pro]{Corollary}

\theoremstyle{definition}
\newtheorem{dfn}[pro]{Definition}
\newtheorem{dfns}[pro]{Definitions}

\theoremstyle{remark}
\newtheorem*{rmk}{Remark}

\newcommand{\pmm}{primitive meridian}

\newcommand{\scc}{simple closed curve}
\newcommand{\s}{\Sigma}

\newcommand{\ie}{{\it i.e.}}
\newcommand{\cf}{{\it cf.}}
\newcommand{\sft}{{swallow follow torus}}

\newcommand{\eg}{{\it e.g.}}

\newcommand{\del}{\partial}
\newcommand{\hh}{Heegaard}
\newcommand{\hhs}{Heegaard surface}

\title{Heegaard genus of the connected sum of {m}-small knots}

\date{February 27, 2006}
\address{Department of Mathematics, Nara Women's University
Kitauoya Nishimachi, Nara 630-8506, Japan}
\address{Department of mathematical Sciences, University of
Arkansas, Fayetteville, AR 72701}
\email{tsuyoshi@cc.nara-wu.ac.jp}
\email{yoav@uark.edu}
\author{Tsuyoshi Kobayashi}
\author{Yo'av Rieck}
\thanks{The first names author was supported by Grant-in-Aid for
scientific research, JSPS grant number 00186751. The second named
author was supported in part JSPS (fellow number P00024) and by
the 21st century COE program ``Constitution for wide-angle
mathematical basis focused on knots" (Osaka City University);
leader: Akio Kawauchi.}

\let\mgp=\marginpar \def\marginpar#1{\mgp{\raggedright\tiny #1}}
\let\lbl=\label
\def\label#1{\lbl{#1}\ifinner\else\marginpar{\ref{#1} #1}\ignorespaces\fi}

\def \marginpar#1{}
\begin{document}

\subjclass{57M99, 57M25}%
\keywords{3-manifold, knots, Heegaard splittings, tunnel number}%

\date{\today}%
% ----------------------------------------------------------------
\begin{abstract}

We prove that if
$K_1 \subset M_1,\dots,K_n \subset M_n$ are m-small knots in
closed orientable 3-manifolds then the Heegaard genus of
$E(\#_{i=1}^n K_i)$ is strictly less than the sum of the Heegaard
genera of the $E(K_i)$ ($i=1,\dots,n$) if and only if there exists
a proper subset $I$ of $\{1,\dots,n\}$ so that $\#_{i \in I} K_i$
admits a primitive meridian.  This generalizes the main result of Morimoto in
\cite{morimoto1}.
\end{abstract}
\maketitle
% ----------------------------------------------------------------
\section{Introduction}
\label{sec:introduction}

Let $K_i$ ($i=1,\dots,n$) be knots in closed orientable
3-manifolds $M_i$, and let $K(=\#_{i=1}^n K_i) \subset M
(=\#_{i=1}^n M_i)$ be their connected sum. Let $X = E(K)$ be the
exterior of $K$ and $X_i$ be the exteriors of $K_i$.  For $n=2$,
it is well-known that $g(X) \leq g(X_1) + g(X_2)$ and by
induction we see that for any $n$, $g(X) \leq \s_{i=1}^n g(X_i)$,
where $g(\cdot)$ denotes the Heegaard genus (definitions are given
in Section~\ref{sec:background}). In this paper we study the
phenomenon $g(X) < \s_{i=1}^n g(X_i)$.  We note that (for $n=2$)
the phenomena $g(X) = g(X_1) + g(X_2)$ and $g(X) \leq g(X_1) +
g(X_2) - k$ (for arbitrary  $k$) can both occur (see
\cite{moriah-rubinstein} and \cite{skuma-morimoto-yokata} for the
former, \cite{morimoto-proc-ams} and \cite{kobayashi-degeneration}
and the latter). An instrumental definition for our work is (see
\cite[Definition~4.6]{moriah-new}):

\begin{dfn}
\label{dfn:pm}
A knot exterior $X$ admits a {\it\pmm} if there exists a minimal
genus Heegaard surface $\s \subset X$ that separates $X$ into a
compression body $C$ and a handlebody $V$ so that there exist a
compressing disk $D \subset V$ and a vertical annulus $A \subset
C$ with $A \cap \del X$ a meridian, and $\del D \subset \s$
intersects $A \cap \s$ transversely in one point.  (For the
definition of vertical annulus, see Remarks~\ref{rmk:properties of
CBs}(5).)
\end{dfn}

\noindent For knots admitting primitive meridians we
have:\footnote{Proposition 1.3 of \cite{morimoto1} is stated in 
terms of the {\it 1-bridge genus.}  However, by \cite[Proposition
2.1]{morimoto1} it is equivalent to this proposition. }

\begin{pro}[\cite{morimoto1}, Propositions 1.3 and 2.1]
\label{pro:primitive meridian implies degeneration}

Let $K_1$, $K_2$ be knots in closed orientable 3-manifolds. Let
$X_1$, $X_2$, and $X$ be the exteriors of the knots $K_1$, $K_2$,
and $K_1 \# K_2$ respectively. If $K_1$ or $K_2$ admits a
primitive meridian, then $g(X) \leq g(X_1) + g(X_2)
-1$.
\end{pro}

\noindent Thus, existence of primitive meridian is a sufficient
condition for: $g(X) < g(X_1) + g(X_2)$. A knot is called {\it
meridionally small} (or {\it m-small}) if there is no essential
surface in $E(K)$ with non-empty boundary so that each boundary
component is a meridian. In \cite[Theorem 1.6]{morimoto1} Morimoto
showed that for m-small knots in $S^3$ existence of primitive
meridian is also necessary for $g(X) < g(X_1) + g(X_2)$, and
conjectured that this holds for any two knots in $S^3$:

\begin{cnj}[Morimoto's Conjecture]
\label{cnj:morimoto's conjecture}

Let $K_1$, $K_2$ be knots in $S^3$.  Let $X_1$, $X_2$, and $X$ be
the exteriors of $K_1$, $K_2$, and $K_1 \# K_2$ respectively. Then
$g(X) < g(X_1) + g(X_2)$ if and only if $X_1$ or $X_2$ admits a
primitive meridian.
\end{cnj}

\noindent In this paper we prove the following generalization of
\cite[Theorem 1.6]{morimoto1}:

\begin{thm}
\label{thm:morimoto's conjecture}

Let $K_1,\dots,K_n$ ($n\geq 2$) be m-small knots in closed
orientable 3-manifolds.  Let $X_1,\dots,X_n$, and $X$ be the
exteriors of $K_1, \dots, K_n$, and $\#_{i=1}^n K_i$ respectively.
Then $g(X) < \s_{i=1}^n g(X_i)$ if and only if there exists $I
\neq \emptyset$ a proper subset of $\{1,\dots,n\}$ such that the
exterior of $\#_{i \in I} K_i$ admits a primitive meridian.
\end{thm}

\begin{rmkks}{\rm
\label{rmk:msy}

\begin{enumerate}
\item One might hope for the following simple statement as the
conclusion of Theorem~\ref{thm:morimoto's conjecture}: given knots
$K_1,\dots,K_n$, $g(X) < \s_{i=1}^n g(X_i)$ holds if and only if
for some $i$, $K_i$ admits a primitive meridian. However, the
following examples show that this is not the case. Let $K_m$ ($m
\in \mathbb Z$) be the knots introduced by
Morimoto--Sakuma--Yokota in \cite{skuma-morimoto-yokata}.  It is
shown in \cite{skuma-morimoto-yokata} that for each $m$,
$g(E(K_m)) = 2$, $K_m$ does not admit a primitive meridian, and
$g(E(K_m \# K_m)) = 4$.   On the other hand, in \cite[Section
5]{kobayashi-rieck-growth-rate} it is shown that $E(K_m \# K_m)$
admits a primitive meridian. Hence, by Proposition
\ref{pro:primitive meridian implies degeneration}, $g(E(K_m \# K_m
\# K_m)) < 4+2 = 3g(E(K_m))$.

\item Many of the ideas used 
for the proof of 
Theorem~\ref{thm:morimoto's
conjecture} can be found in \cite{morimoto2} where Morimoto proved
\cite[Theorem 1.1]{morimoto2} that if $K_i$ ($i=1,\dots,n$) are
m-small knots in closed orientable 3-manifolds $M_i$ so that no
$M_i$ contains a lens space summand, then $g(X) \geq \s_{i=1}^n
g(X_i) - (n-1).$  In this paper we retrieve this result
(Corollary~\ref{cor:bounding-degeneration}) without the
restriction on the summands of the manifolds $M_i$.
\end{enumerate} }
\end{rmkks}

The ``if'' part of Theorem~\ref{thm:morimoto's conjecture} is an easy
consequence of Proposition~\ref{pro:primitive meridian implies degeneration}.
We now sketch the 
proof of the ``only if'' part and introduce the main tools applied.  The proof
is an induction on the number of summands $n$.  Our first tool, the Swallow
Follow Torus Theorem (Theorem~\ref{thm:sft}), implies that if $K_1,\dots,K_n$
($n \geq 2$) are m-small knots then any Heegaard surface for $X(=E(\#_{i=1}^n
K_i))$ weakly reduces to a swallow follow torus (say $T$). In particular,
$X$ does not admit a strongly irreducible Heegaard splitting.  It follows from
the definition of swallow follow torus that $T$ decomposes
$X$ as $X = X_I^{(1)} \cup_T X_J$, with $I \subset \{1,\dots,n\}$ 
a non-empty proper subset, 
$X_I = E(\#_{i \in I} K_i)$, $X_J = E(\#_{i \not\in I} K_i)$, and
$X_I^{(1)}$ is obtained from $X_I$ by drilling out a simple closed curve
parallel to the meridian of $K_I$.  
If $g(X) = g(X_I) + g(X_J)$ then either $g(X_I) < \s_{i \in I} g(E(K_i))$ or
$g(X_J) < \s_{i \in J} g(E(K_i))$ (say the former).  By
induction, Theorem~\ref{thm:morimoto's conjecture} holds for $X_I$ and hence 
for some $I' \subset I$, $E(\#_{i \in I'} K_i)$ admits a primitive meridian,
and we are done.  Thus, we may assume that $g(X) < g(X_I) + g(X_J)$. 
In that case, Corollary~\ref{cor:when can the genus go down} implies that
$g(X_I^{(1)}) = g(X_I)$.

Our second tool, the Hopf-Haken Annulus Theorem
(Theorem~\ref{thm:hopf-annulus}), studies minimal genus
Heegaard splitting of $X_I^{(1)}$.   Denote the boundary of $X_I^{(1)}$ by 
$\partial X_I \cup T$.  Note that $X_I^{(1)}$ admits an essential
annulus (say $A$) connecting $\partial X_I$ to $T$, so that $A \cap
\partial X_I$ is a meridian of $X_I$ and $A \cap T$ is a longitude of
$T$.  The Hopf--Haken Annulus Theorem implies that if $K_I
%%%% Yo'av 2/13/2006 ADDED:
= \#_{i \in I} K_i
%%%%%TO HERE
$ is the
connected sum of m-small knots then there exists a minimal genus
Heegaard surface for $X_I^{(1)}$ (say $\s$) that intersects $A$ in a single
simple closed curve that is essential in both $A$ and the Heegaard
surface.  Applying this to $X_I^{(1)}$ described in
the previous paragraph, since $g(X_I) = g(X_I^{(1)})$ we have that
$\s$ is also a minimal genus Heegaard surface for
$X_I$. It is then easy to see (Proposition~\ref{pro:hopf-haken implies
primitive meridian}) that $\s$ fulfills the conditions of
Definition~\ref{dfn:pm} and hence $X_I$ admits a primitive meridian, proving 
Theorem~\ref{thm:morimoto's conjecture}.
Thus, the Hopf--Haken Annulus Theorem connects the Swallow Follow
Torus Theorem to the existence of a primitive meridian.

\smallskip \noindent
{\bf Acknowledgement:} We would like to thank Kanji Morimoto for
helpful conversations and the referees for many helpful comments.
We thank Sean Bowman for a careful reading of an early manuscript.
The second named author: this research started while I was a JSPS
fellow at Nara Women's University and was finished while I was
member of $21^{\mbox{st}}$-century COE program ``Constitution for
wide-angle mathematical basis focused on knots" at Osaka City
University, leader: Akio Kawauchi. I would like to thank both
universities, the math departments, Tsuyoshi Kobayashi and Akio
Kawauchi for wonderfully warm hospitality.
\section{Background}
\label{sec:background}

Throughout this paper we work in the category of smooth manifolds.
All manifolds considered are assumed to be orientable.  For a
submanifold $H$ of a manifold $M$, $N(H,M)$ denotes a regular
neighborhood of $H$ in $M$.  When $M$ is understood from context
we often abbreviate $N(H,M)$ to $N(H)$.  Let $F$ be a surface.
Then a surface in $F \times [0,1]$ is called {\it vertical} if it
is ambient isotopic to a surface of the form $c \times [0,1]$,
$c\subset F$ an embedded 1-manifold. Let $N$ be a manifold
embedded in a manifold $M$ with $\dim(N) = \dim(M)$.  Then the
frontier of $N$ in $M$ is denoted by $\mbox{Fr}_M(N)$. For
definitions of standard terms of 3-dimensional topology we refer
to \cite{hempel} or \cite{jaco}.

The following lemma is a well-known fact.  Informally, it says
that for manifolds with torus boundary compression
implies compression; for a proof see, for example, \cite[Lemma
2.7]{kobayashi-rieck-local-det}.

\begin{lem}
\label{boudary compressible implies compressible}
Let $N$ be a 3-manifold with a torus boundary component $T$. Let
$S$ be a 2-sided surface properly embedded in $N$ such that $S
\cap T$ consists of essential simple closed curves in $T$. Suppose
that there is a boundary compressing disk $\Delta$ for $S$ such
that $\Delta$ compresses $S$ into $T$, \ie , $\Delta \cap
\partial N = \partial \Delta \cap T$ is an arc, say $\alpha$, and
$\Delta \cap S = \partial \Delta \cap S$ is an essential arc in
$S$, say $\beta$, such that $\alpha \cup \beta = \partial \Delta$.
Then we have either one of the following.

\begin{enumerate}
\item $S$ is compressible. Moreover, if $S$ is separating in $N$,
then the compression occurs in the same side as the boundary
compression. \item $S$ is an annulus; moreover, when $N$ is
irreducible, $S$ is boundary parallel.
\end{enumerate}
\end{lem}

Let $S$ be a surface properly embedded in a 3-manifold $N$.  Then
$S$ is called {\it essential} if $S$ is incompressible and not
boundary parallel.  If $N$ is irreducible and $\partial N$
consists of tori then by Lemma~\ref{boudary compressible implies
compressible}, $S$ being essential implies $S$ is not boundary
compressible.

A 3-manifold $C$ is called a \em compression body \em if there is
a compact, connected, closed, orientable surface $F$ such that $C$
is obtained from $F \times [0,1]$ by attaching 2-handles along
mutually disjoint \scc s in $F \times \{1\}$ and capping off the
2-sphere boundary components using 3-handles.  Then $F \times
\{0\}$ is denoted by $\del_+ C$ and $\del_- C = \del C \setminus
\del_+ C$.  A compression body is called a \em handlebody \em if
$\del_- C = \emptyset$.  A compression body $C$ is called \em
trivial \em if $C$ is homeomorphic to $F \times [0,1]$ with
$\del_+F = F \times \{0\}$ and $\del_-F = F \times \{1\}$.  A disk
$D$ embedded in a compression body $C$ is called a \em meridian
disk \em if $\del D$ is non-trivial in $\del_+ C$.  (We note that
$\del_- C$ is incompressible.) By extending the cores of the
2-handles in the definition of a compression body vertically in $F
\times [0,1]$ we obtain a union of mutually disjoint meridian
disks in $C$, say $\mathcal{D}$, such that the manifold obtained
from $C$ by cutting along $\mathcal{D}$ is homeomorphic to a union
of $\del_- C \times [0,1]$ and a (possibly empty) union of 3-balls
corresponding to the attached 3-handles.  This gives a dual
description of compression bodies, that is, a connected 3-manifold
$C$ is a compression body if there exists a closed (not
necessarily connected) surface $F$ without 2-sphere components and
a (possibly empty) union of 3-balls (say $\mathcal{B}$) so that
$C$ is obtained from $(F \times [0,1]) \cup \mathcal{B}$ by
attaching 1-handles along $(F \times \{0\}) \cup \del
\mathcal{B}$. We may regard each component of $\mathcal{B}$ as a
0-handle. Hence $C$ admits a handle decomposition $(F \times
[0,1]) \cup 0\mbox{-handles } \cup 1\mbox{-handles}$. We note that
$\del_- C$ is the surface corresponding to $F \times \{1\}$.

\begin{rmkks}
\label{rmk:properties of CBs}

{\rm  The following properties are known for compression bodies
(the proofs are omitted or outlined below):}

\begin{enumerate}{\rm
    \item Compression bodies are irreducible.

    \item Let $\mathcal{D}$ be a union of mutually disjoint
    meridian disks of a compression body $C$, and $C'$ a component
    of the manifold obtained by cutting $C$ open along
    $\mathcal{D}$.  Then $C'$ inherits a compression body
    structure from $C$, \ie ,\ $C'$ is a compression body such that
    $\del_- C' = \del_- C \cap C'$ and $\del_+ C' = (C' \cap \del_+ C)
    \cup \mbox{Fr}_C(C')$.

    \item Let $S$ be an incompressible surface in $C$ such that
    $\del S \ne \emptyset$, and
    $\del S \subset \del_+ C$.  Suppose that $S$ is not a union of
    meridian disks.  Then $S$ is boundary compressible into $\del_+ C$
    in $C$, \ie ,\ there exists a disk $\Delta$ in $C$ such that $\Delta
    \cap S = \del \Delta \cap S = a$ is an essential arc in $S$
    and $\Delta \cap \del C = \Delta \cap \del_+ C = \mbox{cl}(\del \Delta
    \setminus
    a)$.

    \item Let $A$ be a surface properly embedded in $C$
    such that $\del A \subset \del_+ C$, and each component of $A$
    is an essential annulus.
    Then by boundary compressing an outermost component
    of $A$ (see (3)
    above) we obtain a meridian disk $D$.  By a tiny isotopy we
    may assume that $D \cap A = \emptyset$.  This implies the
    following: let $A \subset C$ be as above.  Then there is
    a meridian disk $D$ of $C$ such that $D \cap A = \emptyset$.

    \item Let $S$ be an incompressible, boundary incompressible
    surface with non-empty boundary
    properly embedded in $C$.  Then each component of $S$
    is either a meridian disk or an annulus $A$ such that one
    component of $\del A$ is in $\del_+ C$, and the other in
    $\del_- C$.  (Such an incompressible annulus is called \em
    vertical.\em)

    \item Let $S$ be an  incompressible (possibly closed) surface
    properly embedded in
    $C$ such that $S \cap \del_+ C = \emptyset$.
    Then each component of $S$ is parallel to a subsurface of  $\del_- C$.}
\end{enumerate}
\end{rmkks}

Let $N$ be a 3-manifold and $F_1$, $F_2$ a partition of $\partial
N$ (possibly, $F_1 = \emptyset$ or $F_2 = \emptyset$).

\begin{dfn}
\label{dfn:heegaard splitting}

We say that a decomposition $C_1 \cup_\s C_2$ is a {\it Heegaard
splitting} of $(N;F_1,F_2)$ if the following conditions hold:

\begin{enumerate}
    \item $C_i$ ($i=1,2$) is a compression body,
    \item $C_1 \cup C_2 = N$,
    \item $C_1 \cap C_2 = \del_+ C_1 = \del_+ C_2 = \s$, and
    \item $\del_- C_i = F_i$ ($i=1,2$).
\end{enumerate}
\end{dfn}

\noindent The surface $\s$ is called a \em Heegaard surface \em
of $(N;F_1,F_2)$.  If the boundary partition of $N$ is irrelevant,
we omit condition (4) in Definition~\ref{dfn:heegaard splitting}
and say that $\s$ is a Heegaard surface of $N$. The genus of the
Heegaard splitting is the genus of $\s$, denoted by $g(\s)$ or
$g(C_1 \cup_\s C_2)$. The minimal genus of all Heegaard splittings
of $(N;F_1,F_2)$ is called the {\it(Heegaard) genus } of
$(N;F_1,F_2)$ and is denoted by $g(N;F_1,F_2)$. The minimal genus
of all Heegaard splittings of $N$ is called the {\it (Heegaard)
genus } of $N$ and is denoted by $g(N)$.

\begin{dfns}
\label{dfn:reducibility of heegaard splittings}

\begin{enumerate}
    \item A Heegaard splitting $C_1 \cup_\s C_2$ is called
    \em reducible \em if there exist meridian disks $D_1 \subset
    C_1$ and $D_2 \subset C_2$ such that $\del D_1 = \del D_2$.
    Otherwise, the Heegaard splitting is called \em irreducible. \em
    \item A Heegaard splitting $C_1 \cup_\s C_2$ is called
    \em weakly reducible \em if there exist meridian disks $D_1 \subset
    C_1$ and $D_2 \subset C_2$ such that $\del D_1 \cap \del D_2 =
    \emptyset$.  Otherwise the splitting is called \em strongly
    irreducible.\em
    \item A Heegaard splitting $C_1 \cup_\s C_2$ is called
    \em stabilized \em if there exist meridian disks $D_1 \subset
    C_1$ and $D_2 \subset C_2$ such that $\del D_1$ intersects
    $\del D_2$ transversely in one point.
    Otherwise the splitting is said to be \em non stabilized.\em
    \item A Heegaard splitting $C_1 \cup_\s C_2$ is called
    \em trivial \em if $C_1$ or $C_2$ is a trivial compression body.
\end{enumerate}
\end{dfns}

Let $C_1 \cup_\s C_2$ be a Heegaard splitting of $(N;F_1,F_2)$.
Recall that $C_1$ was obtained from $(F_1 \times [0,1]) \cup$
0-handles by attaching 1-handles, and that $C_2$ was obtained from
$\del_+ C_2 \times [0,1]$ by attaching 2- and 3-handles.  Then, by
using an isotopy which pushes $\del_+ C_2 \times [0,1]$ out of
$C_2$, we identify $\del_+ C_2 \times [0,1]$ with $N(\del_+
C_1,C_1)$. This identification, together with the above handles
gives the following handle decomposition of $N$:
$$N = (F_1 \times [0,1]) \cup 0\mbox{-handles } \cup 1
\mbox{-handles } \cup 2\mbox{-handles } \cup 3\mbox{-handles. }$$
\noindent We say that this handle decomposition is {\it induced} by
$C_1 \cup_\s C_2$.  Suppose that there exists a proper subset of
the 0- and 1-handles so that some subset of the 2- and 3-handles
does not intersect the subset of 0- and 1-handles. Suppose further
that we can rearrange the order of the handles non-trivially to
obtain an increasing sequence of (not necessarily connected)
submanifolds $N_1,\dots ,N_n$ of $N$ such that the following
holds:

$N_0=\emptyset$,

$N_j = N_{j-1} \cup (F_1^{(j-1)} \times [0,1])
\cup \mbox{ 0-handles }
\cup \mbox{ 1-handles } \cup$

$\mbox{\hspace{1in}}\cup \mbox{ 2-handles } \cup
\mbox{ 3-handles} \hspace{3mm} (1 \leq j \leq n),$

\noindent where $\{F_1^{(0)},\dots,F_1^{(n-1)}\}$ is a partition
of the components of $F_1$ (with $F_1^{(j)}$ possibly empty, $0
\leq j \leq n-1$) and the handles appearing in the above come from
a handle decomposition of $N$ induced by $C_1 \cup_\s C_2$, where
each collection of 0- and 3-handles may be empty but the
collections of 1- and 2-handles are never empty.  Suppose further
that this handle decomposition satisfies the following three
conditions:

\begin{enumerate}
    \item $N_1$ is connected.
    \item At each stage $j$ ($1 \leq j \leq n$), let $\del^-_j$
    denote the union of the components of $\del N_{j-1}$ to which
    the 1-handles are attached.  Then $\del_j^- \cup (F_1^{(j-1)}
    \times [0,1]) \cup $ 0-handles $\cup$ 1-handles $\cup$
    2-handles $\cup$ 3-handles is connected, where these handles are the
    handles that appear in the description of $N_j$.
    \item No component of $\del N_j$ is a 2-sphere
    ($j=1,\dots,n-1$).
\end{enumerate}

\noindent We note that if $C_1 \cup_\s C_2$ is irreducible then
condition (3) holds automatically 
(this is a consequence of  Lemma~\ref{lem:weak
reduction that yeilds a 2-sphere} and Proposition~\ref{pro:untel
to a conn sep surfce implies weak reduction} of this paper). 
Then for each $j,$
let $I_j = \del_j^- \times [0,1]$, and $R_j = I_j \cup (F_1^{(j)}
\times [0,1]) \cup \mbox{ 0-handles } \cup \mbox{ 1-handles } \cup
\mbox{ 2-handles } \cup \mbox{ 3-handles }$ (handles that appear
in the description of $N_j$), where $\del_j^-$ $(\subset \del
N_{j-1})$ is identified with $\del_j^- \times \{ 0 \} (\subset
I_j)$. By the above conditions, we see that this handle
decomposition induces a Heegaard splitting, say $C_1^{(j)} \cup
C_2^{(j)}$, of $R_j$.

It is clear that $N$ can be regarded as obtained from
$R_1,\dots,R_n$ by identifying their boundaries. Hence we have
obtained a decomposition of $N$ into submanifolds $R_1,\dots,R_n$
and Heegaard splittings for those submanifolds as follows:
$$(C_1^{(1)}\cup C_2^{(1)})\cup(C_1^{(2)} \cup C_2^{(2)})\cup
                    \cdots \cup (C_1^{(n)} \cup C_2^{(n)}).$$
This decomposition is called an \em untelescoping \em of $C_1
\cup_\s C_2$, originally defined in a slightly different way by Scharlemann
and Thompson in
\cite{scharl-abby}. An untelescoping is called a \em
Scharlemann--Thompson untelescoping \em (or S--T untelescoping) of
the Heegaard splitting $C_1 \cup_\s C_2$ if for each $j$,
$C_1^{(j)} \cup C_2^{(j)}$ is strongly irreducible.

\begin{rmkk}
\label{rmk:scharlemann's ST-untel is different}

{\rm It is well known that every irreducible Heegaard splitting of
a 3-manifold $M$ with incompressible boundary (other than surface
cross interval) admits a S--T untelescoping such that each
Heegaard splitting is non-trivial (see \cite{scharl-abby} and
\cite{scharlemann-handbook}).  In this case, every component of
$\del R_j$ is incompressible in $M$.  If in addition $C_1 \cup_\s
C_2$ is minimal genus then by \cite{eric-agt} each component of
$\del R_j \setminus \del M$ is an essential surface.

}
\end{rmkk}

\smallskip \noindent
{\bf Amalgamation of Heegaard splittings.}

Let $M^{(1)}$ and $M^{(2)}$ be compact, orientable 3-manifolds
with boundary, together with partitions of the boundary components:

$$\del M^{(1)} = F^{(1)} \cup F_1^{(1)} \cup F_2^{(1)},$$

$$\del M^{(2)} = F^{(2)} \cup F_1^{(2)} \cup F_2^{(2)},$$

\noindent
where $F^{(1)}$ and $F^{(2)}$ are non-empty and mutually
homeomorphic.  Let $M$ be a manifold obtained from $M^{(1)}$ and
$M^{(2)}$ by identifying $F^{(1)}$ and $F^{(2)}$.  Then $F$ denotes the
image of $F^{(1)}$ (= the image of $F^{(2)}$) in $M$.

Suppose that we are given Heegaard splittings

$$C_1^{(1)} \cup C_2^{(1)} \mbox{ of }
(M^{(1)};F_1^{(1)},F^{(1)} \cup F_2^{(1)})$$

\noindent
and---

$$
C_1^{(2)} \cup C_2^{(2)}\mbox{ of } (M^{(2)};F^{(2)} \cup
F_1^{(2)},F_2^{(2)}).$$

\noindent
Recall that there are handle decompositions:

$$M^{(1)} = (F_1^{(1)} \times [0,1]) \cup \mathcal{H}_0^{(1)}
\cup \mathcal{H}_1^{(1)}\cup \mathcal{H}_2^{(1)}\cup
\mathcal{H}_3^{(1)},$$

\noindent
and---

$$M^{(2)} = (F^{(2)} \times [0,1]) \cup (F_1^{(2)} \times [0,1])
\cup \mathcal{H}_0^{(2)} \cup \mathcal{H}_1^{(2)}\cup
\mathcal{H}_2^{(2)}\cup \mathcal{H}_3^{(2)},$$

\noindent where the $\mathcal{H}_i^{(j)}$ denotes the handles of
index $i$ induced by the Heegaard splitting $C_1^{(j)} \cup
C_2^{(j)}$. Hence we have:

$$M = (F_1^{(1)} \times [0,1]) \cup \mathcal{H}_0^{(1)}
\cup \mathcal{H}_1^{(1)}\cup \mathcal{H}_2^{(1)}\cup
\mathcal{H}_3^{(1)}$$

$$\cup (F^{(2)} \times [0,1]) \cup (F_1^{(2)} \times [0,1])
\cup \mathcal{H}_0^{(2)} \cup
\mathcal{H}_1^{(2)}\cup \mathcal{H}_2^{(2)} \cup
\mathcal{H}_3^{(2)},$$

\noindent Then, via ambient isotopy, we push $M^{(2)}$ slightly
into $M^{(1)}$ so that $F^{(2)} \times [0,1]$ is identified with
$N(F^{(1)}, M^{(1)})$. Hence $F^{(1)}$ ($\subset \del M^{(1)}$) is
regarded as $F^{(2)} \times \{0\} \subset \mbox{int} M^{(2)}$.
(Note that in the handle decomposition $M^{(2)} = (F^{(2)} \times
[0,1]) \cup (F_1^{(2)} \times [0,1]) \cup \mathcal{H}_0^{(2)} \cup
\mathcal{H}_1^{(2)}\cup \mathcal{H}_2^{(2)}\cup
\mathcal{H}_3^{(2)}$, the handles of $\mathcal{H}_1^{(2)}$ are
attached to $(F^{(2)} \times \{0\}) \cup (F_1^{(2)} \times \{0\})
\cup \del \mathcal{H}_0^{(2)}$.) Here, via isotopy, we may suppose
that $\mathcal{H}_1^{(2)}$ is disjoint from $\mathcal{H}_1^{(1)}
\cup \mathcal{H}_2^{(1)}$ and that $\mathcal{H}_2^{(2)}$ is
disjoint from $\mathcal{H}_2^{(1)}$. This implies that we can
change the order of the handles so that:

$$M = (F_1^{(1)} \times [0,1]) \cup (F_1^{(2)} \times [0,1])
\cup(\mathcal{H}_0^{(1)} \cup \mathcal{H}_0^{(2)})$$

$$\cup(\mathcal{H}_1^{(1)}\cup\mathcal{H}_1^{(2)})
\cup(\mathcal{H}_2^{(1)}\cup \mathcal{H}_2^{(2)} ) \cup
(\mathcal{H}_3^{(1)} \cup \mathcal{H}_3^{(2)}).$$

By using this handle decomposition we can obtain a Heegaard
splitting, say $C_1 \cup C_2$, of $M$.  We say that $C_1 \cup C_2$
is obtained from $C_1^{(1)} \cup C_2^{(1)}$ and $C_1^{(2)} \cup
C_2^{(2)}$ by an \em amalgamation, \em or that $C_1 \cup C_2$ is
an amalgamation of $C_1^{(1)} \cup C_2^{(1)}$ and $C_1^{(2)} \cup
C_2^{(2)}$.  We note that this definition of amalgamation has
different appearance than that in \cite{schultens-FXS1}.  However,
it is easy to see that they yield the same surface.

By the definition of untelescoping we immediately have the
following:

\begin{pro}

Let $(C_1^{(1)} \cup C_2^{(1)}) \cup \dots \cup (C_1^{(n)} \cup
C_2^{(n)})$ be an untelescoping of $C_1 \cup C_2$.  Then $C_1 \cup
C_2$ is obtained from $C_1^{(1)} \cup C_2^{(1)}, \dots, C_1^{(n)}
\cup C_2^{(n)}$ by a sequence of amalgamations.
\end{pro}

Proof of the following lemma (\eg\ Remark 2.7 of
\cite{schultens-FXS1} for the case $m=1$) is elementary and is
left to the reader.

\begin{lem}
\label{lem:genus after amalgamation}

Let $C_1 \cup C_2$, $C_1^{(1)} \cup C_2^{(1)}$ and $C_1^{(2)} \cup
C_2^{(2)}$ be as above.  Let $F_1,\dots,F_m$ be the components of
$F \subset M$. Then
$$g(C_1 \cup C_2) = g(C_1^{(1)} \cup C_2^{(1)}) + g(C_1^{(2)}
\cup C_2^{(2)}) - \s_{i=1}^m g(F_i) + (m-1).$$

\end{lem}

\begin{pro}
\label{pro:amalgamation of minimal genus}

Let $M$ be a compact 3-manifold with partition of the boundary
components:
$\del M = (F_1^{(1)} \cup F_1^{(2)}) \cup (F_2^{(1)} \cup F_2^{(2)}).$
Suppose that there exists a minimal genus Heegaard splitting $C_1
\cup C_2$ of $(M;F_1^{(1)} \cup F_1^{(2)}, F_2^{(1)} \cup
F_2^{(2)})$ such that $C_1 \cup C_2$ admits an untelescoping

$$(C_1^{(1)} \cup C_2^{(1)}) \cup_F (C_1^{(2)} \cup C_2^{(2)})$$

\noindent
with the following properties:

\begin{itemize}
    \item $\del_-C_1^{(1)} = F_1^{(1)}$ and $\del_-C_2^{(2)} =
    F_2^{(2)}$, and---
    \item $\del_-C_2^{(1)} = F^{(1)} \cup F_2^{(1)}$ and $\del_-C_1^{(2)} =
    F^{(2)} \cup F_1^{(2)}$,
\end{itemize}

\noindent
where $F^{(1)}$ and $F^{(2)}$ are the surfaces that are identified
in $M$ as $F$.
Then $C_1^{(1)} \cup C_2^{(1)}$ and $C_1^{(2)} \cup C_2^{(2)}$ are
minimal genus Heegaard splittings of $(M^{(1)};F_1^{(1)},F^{(1)}
\cup F_2^{(1)})$ and  $(M^{(2)};F^{(2)} \cup F_1^{(2)},F_2^{(2)})$
respectively.  In particular, we have the following identity:

\begin{eqnarray*}
  g(M;F_1^{(1)} \cup F_1^{(2)},F_2^{(1)} \cup F_2^{(2)})
        &=& g(M^{(1)};F_1^{(1)},F_2^{(1)} \cup F^{(1)}) \\
  && + g(M^{(2)};F_1^{(2)}\cup F^{(2)} ,F_2^{(2)})\\
   && - \s_{i=1}^m g(F_i) + (m-1),
\end{eqnarray*}

\noindent
where $F_1,\dots,F_m$ are the components of $F$.
\end{pro}

\begin{proof}
If either $C_1^{(1)} \cup C_2^{(1)}$ or $C_1^{(2)} \cup C_2^{(2)}$
is not a minimal genus Heegaard splitting of
$(M^{(1)};F_1^{(1)},F_2^{(1)} \cup F^{(1)})$ or
$(M^{(2)};F_1^{(2)}\cup F^{(2)},F_2^{(2)})$ (respectively), then
by amalgamating minimal genus Heegaard splittings of
$(M^{(1)};F_1^{(1)},F^{(1)} \cup F_2^{(1)})$ and
$(M^{(2)};F_1^{(2)}\cup F^{(2)},F_2^{(2)})$ we see, by
Lemma~\ref{lem:genus after amalgamation}, that we obtain a
Heegaard splitting of $(M;F_1^{(1)} \cup F_1^{(2)},F_2^{(1)} \cup
F_2^{(2)})$ with genus lower than $g(C_1 \cup C_2)$, contradicting
the fact that the $C_1 \cup C_2$ is a minimal genus Heegaard
splitting of $(M;F_1^{(1)} \cup F_1^{(2)},F_2^{(1)} \cup
F_2^{(2)})$. Hence $C_1^{(1)} \cup C_2^{(1)}$ and $C_1^{(2)} \cup
C_2^{(2)}$ are minimal genus Heegaard splittings of
$(M^{(1)};F_1^{(1)},F^{(1)} \cup F_2^{(1)})$ and
$(M^{(2)};F_1^{(2)}\cup F^{(2)},F_2^{(2)})$ respectively. This
together with Lemma~\ref{lem:genus after amalgamation} implies:

\begin{eqnarray*}
  g(M;F_1^{(1)} \cup F_1^{(2)},F_2^{(1)} \cup F_2^{(2)})
        &=& g(M^{(1)};F_1^{(1)},F_2^{(1)} \cup F^{(1)}) \\
  && + g(M^{(2)};F_1^{(2)}\cup F^{(2)} ,F_2^{(2)})\\
   && - \s_{i=1}^m g(F_i) + (m-1).
\end{eqnarray*}
\end{proof}

Recall that $g(M)$ denotes the genus of the minimal genus Heegaard
splitting of $M$ for all possible partitions of the components of
$\del M$.

\begin{pro}
\label{prop:amalgamation of minimal genus: connected case}

Let $N$ be a compact orientable 3-manifold with a minimal genus
Heegaard splitting $C_1 \cup_\s C_2$ admitting an untelescoping
$(C_1^{(1)} \cup C_2^{(1)}) \cup_F (C_1^{(2)} \cup C_2^{(2)})$,
where $F$ is a connected surface.  Let $M^{(1)} = C_1^{(1)} \cup
C_2^{(1)}$ and $M^{(2)} = C_1^{(2)} \cup C_2^{(2)}$.  Then
$C_1^{(1)} \cup C_2^{(1)}$ and $C_1^{(2)} \cup C_2^{(2)}$ are
minimal genus Heegaard splittings of $M^{(1)}$ and $M^{(2)}$
respectively.  In particular we have the following equality:

$$g(M) = g(M^{(1)}) + g(M^{(2)}) - g(F).$$
\end{pro}

\begin{proof}

Suppose that either $C_1^{(1)} \cup C_2^{(1)}$ or $C_1^{(2)} \cup
C_2^{(2)}$, say $C_1^{(1)} \cup C_2^{(1)}$, is not a minimal genus
Heegaard splitting.  Let $\bar C_1^{(1)} \cup \bar C_2^{(1)}$ be a
minimal genus Heegaard splitting of $M^{(1)}$. Since $F$ is
connected, by exchanging $\bar C_1^{(1)}$ and $\bar C_2^{(1)}$ if
necessary, we may assume that $F \subset \del_- \bar C_2^{(1)}$.
Hence we can amalgamate  $\bar C_1^{(1)} \cup \bar C_2^{(1)}$ and
$C_1^{(2)} \cup C_2^{(2)}$ to obtain a Heegaard splitting, say
$\bar C_1 \cup_{\bar\s} \bar C_2$ of $M$.  Here we note that
$g(\bar\s) < g(\s)$ by Lemma~\ref{lem:genus after amalgamation}, a
contradiction. Hence $C_1^{(1)} \cup C_2^{(1)}$ and $C_1^{(2)}
\cup C_2^{(2)}$ are minimal genus Heegaard splittings. This
together with Lemma~\ref{lem:genus after amalgamation} implies:

$$ g(M) = g(M^{(1)}) + g(M^{(2)}) - g(F).$$
\end{proof}

\noindent{\bf Weak reduction of Heegaard splittings.}

Let $M$ be a compact, orientable 3-manifold, and $C_1 \cup_\s C_2$
a Heegaard splitting of $M$.  Let $\Delta = \Delta_1 \cup
\Delta_2$ be a weakly reducing collection of disks for $\s,$ \ie
,\ $\Delta_i$ ($i=1,2$) is a union of mutually disjoint, non-empty
meridian disks for $C_i$ such that $\Delta_1 \cap \Delta_2 =
\emptyset$.  Then $\s(\Delta)$ denotes the surface obtained from
$\s$ by compressing along $\Delta$, \ie, $\s(\Delta) = (\s
\setminus (N(\Delta_1, C_1) \cup N(\Delta_2, C_2))) \cup
\mbox{Fr}_{C_1}N(\Delta_1, C_1) \cup  \mbox{Fr}_{C_2} N(\Delta_2,
C_2)$. Let $\widehat{\s}(\Delta)$ be the surface obtained from
${\s}(\Delta)$ by removing all the components that are contained
in $C_1$ or $C_2$. We call $\widehat{\s}(\Delta)$ \em the surface
obtained from $\s$ by weakly reducing along $\Delta$ \em (or \em
the surface obtained from $\s$ by a weak reduction along
$\Delta$\em). The next lemma is well known (for a proof see, for
example, \cite[Lemma 4.1]{kobayashi-2-bridge}).

\begin{lem}
\label{lem:weak reduction that yeilds a 2-sphere}

If there is a 2-sphere component of $\widehat{\s}(\Delta)$, then
$C_1 \cup_\s C_2$ is reducible.
\end{lem}

Suppose that no component of  $\widehat{\s}(\Delta)$ is a
2-sphere.  (By Lemma~\ref{lem:weak reduction that yeilds a
2-sphere} we see that when $C_1 \cup_\s C_2$ is irreducible this
is always satisfied.)  The following argument can be found in
\cite{boilau-otal} (see also \cite[Section 4]{kobayashi-2-bridge}).

Let $M_1,\dots,M_n$ be the closures of the components of $M
\setminus \widehat{\s}(\Delta)$, and let $M_{i,j} = M_i \cap C_j$,
($i=1,\dots,n$, $j=1,2$).  Then we can show (for a proof see, for example,
\cite[Lemma 4.3]{kobayashi-2-bridge}):

\begin{pro}
\label{kobayashi-2-bridge4.3}

For each $j$ one of the following holds:

\begin{enumerate}
    \item $M_{j,2} \cap \s \subset \mbox{\rm int}(M_{j,1} \cap \s)$ and
    $M_{j,1}$ is connected,
    \item $M_{j,1} \cap \s \subset \mbox{\rm int}(M_{j,2} \cap \s)$ and
    $M_{j,2}$ is connected.
\end{enumerate}
\end{pro}

Suppose that $M_j$ satisfies the conclusion (1) ((2) resp.) of
Proposition~\ref{kobayashi-2-bridge4.3}. 
Note that $M_{j,1}$ ($M_{j,2}$ rersp.) is a compression body
((2) of Remarks~\ref{rmk:properties of CBs}). 
Let $C_{j,1} =
\mbox{cl}(M_{j,1} \setminus N(\del_+ M_{j,1},M_{j,1}))$ ($C_{j,2}
= \mbox{cl}(M_{j,2} \setminus N(\del_+ M_{j,2},M_{j,2}))$ resp.),
and $C_{j,2} = N(\del_+ M_{j,1}, M_{j,1}) \cup M_{j,2}$ ($C_{j,1}
= N(\del_+ M_{j,2}, M_{j,2}) \cup M_{j,1}$ resp.). Note that
$C_{j,1}$ ($C_{j,2}$ resp.) is a compression body ((2) of
Remarks~\ref{rmk:properties of CBs}), and that $C_{j,2}$
($C_{j,1}$ resp.) is also a compression body by the definition of
compression body. By using this it is easy to see:

\begin{pro}
\label{pro:weak reduction gives untel}

Each $C_{i,j}$ is a compression body such that $C_{j,1} \cup
C_{j,2}$ gives a Heegaard splitting of $M_j$ ($j=1,\dots ,n$) and
that

$$(C_{1,1} \cup C_{1,2})\cup \cdots \cup (C_{n,1} \cup C_{n,2})$$

\noindent
gives an untelescoping of $C_1 \cup C_2$.
\end{pro}

\noindent{\bf Untelescoping and weak reduction.}

Let $C_1 \cup_\s C_2$ be a Heegaard splitting of $M$ and
$(C_1^{(1)} \cup_{\s^{(1)}} C_2^{(1)}) \cup \cdots \cup (C_1^{(n)}
\cup_{\s^{(n)}} C_2^{(n)})$ an untelescoping of $C_1 \cup_\s
C_2$.  Let $R^{(i)} = C_1^{(i)} \cup C_2^{(i)}$.

\begin{pro}
\label{pro:untel to a conn sep surfce implies weak reduction}
Suppose that there exists an $R^{(k)}$ such that a component of
$\mbox{\rm Fr}_M (R^{(k)})$, say $F$, is separating in $M$. Then
there exists a weak reducing collection of disks for $\s$, say
$\Delta = \Delta_1 \cup \Delta_2$, so that $\widehat{\s}(\Delta)$
is isotopic to $F$.
\end{pro}

\begin{proof}

By exchanging superscripts if necessary, we may assume that $F$
separates $M$ into $R^{(1)} \cup \cdots \cup R^{(k)}$, say
$M^{(1)}$, and $R^{(k+1)} \cup \cdots \cup R^{(n)}$, say
$M^{(2)}$.  Let $\mathcal{C}_1^{(1)} \cup \mathcal{C}_2^{(1)}$
($\mathcal{C}_1^{(2)} \cup \mathcal{C}_2^{(2)}$ resp.) be a
Heegaard splitting of $M^{(1)}$ ($M^{(2)}$ resp.) obtained from
$(C_1^{(1)} \cup C_2^{(1)}),\dots ,(C_1^{(k)} \cup C_2^{(k)})$
($(C_1^{(k+1)} \cup C_2^{(k+1)}),\dots ,(C_1^{(n)} \cup
C_2^{(n)})$ resp.) by a sequence of amalgamations. Then let:

$$(\del_- C_1^{(1)} \times [0,1]) \cup \mathcal{H}_0^{(1)}
\cup \mathcal{H}_1^{(1)} \cup \mathcal{H}_2^{(1)} \cup
\mathcal{H}_3^{(1)},$$

$$(\del_- C_1^{(2)} \times [0,1]) \cup \mathcal{H}_0^{(2)}
\cup \mathcal{H}_1^{(2)} \cup \mathcal{H}_2^{(2)} \cup
\mathcal{H}_3^{(2)}$$

\noindent
be handle decompositions of $M^{(1)}$ and $M^{(2)}$ respectively
that can be used to amalgamate $\mathcal{C}_1^{(1)} \cup
\mathcal{C}_2^{(1)}$ and $\mathcal{C}_1^{(2)} \cup
\mathcal{C}_2^{(2)}$, \ie, $\mathcal{H}_1^{(2)} \cap
(\mathcal{H}_1^{(1)} \cup \mathcal{H}_2^{(1)}) = \emptyset$ and
$\mathcal{H}_2^{(2)} \cap \mathcal{H}_2^{(1)} = \emptyset$. Hence
the handle decomposition:

$$(\del_-C_1^{(1)}) \cup ((\del_-C_1^{(2)} \setminus F) \times
[0,1]) \cup (\mathcal{H}_0^{(1)} \cup \mathcal{H}_0^{(2)}) \cup$$

\begin{flushright}
$(\mathcal{H}_1^{(1)} \cup \mathcal{H}_1^{(2)}) \cup
(\mathcal{H}_2^{(1)} \cup \mathcal{H}_2^{(2)}) \cup
(\mathcal{H}_3^{(1)} \cup \mathcal{H}_3^{(2)})$
\end{flushright}
corresponds to $C_1 \cup C_2$.  Then take a weak reducing
collection $\Delta = \Delta_1 \cup \Delta_2$ of $C_1 \cup C_2$
such that $\Delta_1$ corresponds to the co-cores of
$\mathcal{H}_1^{(2)}$ and $\Delta_2$ corresponds to the cores of
$\mathcal{H}_2^{(1)}$.  It is easy to see that the weak reduction of
$\s$ along $\Delta$ yields a surface isotopic to $F$, thus proving
Proposition~\ref{pro:untel to a conn sep surfce implies weak reduction}.
\end{proof}

\noindent{\bf Minimal genus Heegaard splittings of (disk with
$n$-holes)$\times S^1$.}

Let $D(n)$ be the 3-manifold (disk with $n$ holes)$\times S^1$
with $n \ge 2$\label{d(n)}.
Let $T_0,T_1,\dots,T_{n}$ be the boundary components of
$D(n)$.  For each pair of integers $p,q (\geq 0)$ with $p+q = n+1$
denote by $g_{p,q}$ the Heegaard genus of $D(n)$ partitioning the
boundary into $p$ and $q$ components, that is:

$$g_{p,q} = g(D(n);\cup_{i=0}^{p-1} T_i, \cup_{i=p}^{n}
T_i).$$

The next proposition is required in Section~\ref{sec:hopf-annulus}
and is easily deduced from \cite{schultens-SFS}.

\begin{pro}
\label{pro:genus of D(n)}
$g_{0,n+1}=g_{n+1,0}= n+1$, otherwise $g_{p,q}= n$.
\end{pro}

Moreover, by \cite{schultens-SFS} we see that minimal genus
Heegaard splittings are given as follows: let $Q$ be a disk with
$n$ holes in $D(n)$, which is a cross section of (disk with $n$
holes)$\times S^1$. Let $\alpha_i$ and $\beta_i$ be arcs properly
embedded in $Q$ as in Figure~2.1. Let $\del_0,\del_1,\dots,\del_n$
be boundary components of $Q$ such that $T_i = \del_i \times S^1$.

\begin{figure}[ht]
\begin{center}
\includegraphics[width=6cm, clip]{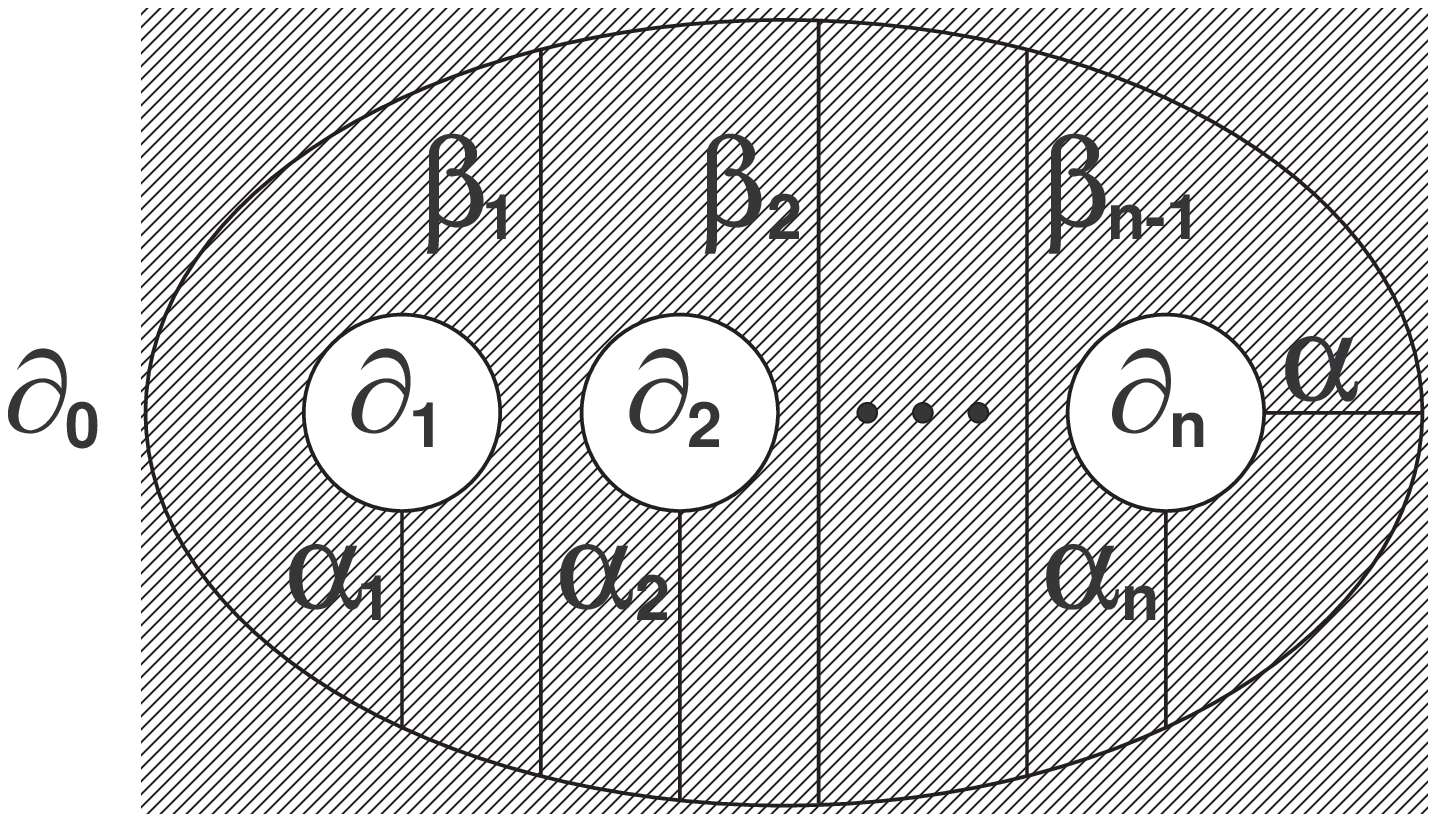}
\end{center}
\begin{center}
Figure 2.1.
\end{center}
\end{figure}

Suppose that $\{p,q\}=\{0,n+1\}$.  Let $C_1 = N(\del D(n) \cup (\cup_{i=1}^{n}
\alpha_i))$, and let $C_2 = \mbox{cl}(D(n) \setminus C_1)$.  Then $C_1 \cup
C_2$ gives a minimal genus Heegaard splitting of $(D(n);\del
D(n),\emptyset)$.

Suppose that $\{p,q\} \neq \{0,n+1\}$.  Let
$C_1=N((\cup_{i=0}^{p-1}T_i) \cup (\cup_{i=1}^{p-1} \alpha_i) \cup
(\cup_{i=p}^{n-1} \beta_i))$, and $C_2 = \mbox{cl}(D(n) \setminus
C_1)$.  Then $C_1 \cup C_2$ gives a minimal genus Heegaard
splitting of $(D(n);\cup_{i=0}^{p-1} T_i,\cup_{i=p}^n T_i)$.

\section{Haken Annuli}
\label{sec:haken annuli-definitions}

In this section we define Haken annuli and prove the basic facts
about their behavior under amalgamation.

\begin{dfn}
\label{dfn:haken annulus}

Let $C_1 \cup_\s C_2$ be a Heegaard splitting of $M$. An essential
annulus $A \subset M$ is called a \em Haken annulus \em for $C_1
\cup_\s C_2$ if $A \cap \s$ consists of a single \scc\ that is
essential in $A$.

\end{dfn}

\begin{rmkks}
\label{rmk:D(n) admits a haken annulus} {\rm
\begin{enumerate}
\item Note the similarity of Haken annulus to the Haken sphere
  \cite{haken} (see also \cite[Chapter 5]{jaco}) and Haken disk
  \cite{casson-gordon}.
\item Let $D(n)$ be as in the end of Section~\ref{sec:background}.
 Observe that if $\{p,q\} \neq
  \{0,n+1\}$ then there is a minimal genus Heegaard surface of $(D(n);
  \cup_{i=0}^{p-1} T_i, \cup_{i=p}^n T_i)$ with a Haken annulus of the
  form $\alpha \times S^1$, where $\alpha$ is an essential arc as in
  Figure~2.1.
\end{enumerate}    }
\end{rmkks}

Given $C_1 \cup_\s C_2$ we described in Section~\ref{sec:background}
how to obtain a (not necessarily unique) handle
decomposition for $M$ of the form:

$$ M = (\del_- C_1 \times [0,1]) \cup \mathcal{H}^0 \cup \mathcal{H}^1
  \cup \mathcal{H}^2 \cup \mathcal{H}^3,$$

\noindent
where $\mathcal{H}^i$ are handles of index $i$.

We say that an annulus $A$ properly embedded in $M$ is 
{\it vertical} with respect to $C_1 \cup_\s C_2$ 
if $A \cap (\cup_{i=1}^3  \mathcal{H}^i) = \emptyset$
(hence $A \subset \partial_-C \times [0,1]$), $A$ is vertical in
$\del_-C \times [0,1]$, and $A \cap (\del_-C_1 \times \{0\})$ is
an essential \scc\ in $\del_-C_1 \times \{0\}$.

From the definition of Haken annulus and Heegaard splitting
induced by handle decomposition above, it is easy to see that:

\begin{lem}
\label{lem:obvious Haken annulus}
If $A$ is vertical 
with respect to $C_1 \cup_\s C_2$ 
then $A$ is a Haken annulus for the Heegaard splitting
$C_1 \cup_\s C_2$.
\end{lem}

Conversely, we have the following:

\begin{lem}
\label{lem:haken annulus disjoint from handles}

Suppose that there exists a Haken annulus $A$ for $C_1 \cup_\s
C_2$.  Then there exists a handle decomposition (with no 0- or
3-handles):

$$ M = (\del_- C_1 \times [0,1]) \cup \mathcal{H}^1
  \cup \mathcal{H}^2 $$

\noindent induced by $C_1 \cup_\s C_2$ such that $A \cap
(\mathcal{H}^1 \cup \mathcal{H}^2) = \emptyset.$
\end{lem}

\begin{proof}

Since $C_1$ is a compression body, there exists a disjoint
collection of meridian disks for $C_1$ (say $\mathcal{D}_1$) such
that $C_1$ cut open along $\mathcal{D}_1$ is $\del_- C_1 \times
[0,1]$ where each component of $\mathcal{D}_1$ can be regarded as
the co-core of a 1-handle for $C_1$.  (Since $\del_- C_1 \neq
\emptyset$, we may assume $C_1$ has no 0-handles.)

Note that $A \cap C_1$ is a vertical annulus in $C_1$. Then by
applying a standard innermost loop, outermost arc arguments to
$\mathcal{D}_1$ and $A \cap C_1$ we can replace the collection
$\mathcal{D}_1$ (without renaming it) so that $(A \cap C_1) \cap
\mathcal{D}_1 = \emptyset$.  By using the same argument we see
that there exists a set of compressing disks for $C_2$ (say
$\mathcal{D}_2$) such that $\mathcal{D}_2 \cap (A \cap C_2) =
\emptyset$ and $\mathcal{D}_2$ cuts $C_2$ into $\del_- C_2 \times
[0,1]$ (since there exists a Haken annulus, $\del_- C_2 \neq
\emptyset$, hence we may assume no component of $C_2$ cut along
$\mathcal{D}_2$ is a 3-ball.) Each component of $\mathcal{D}_2$
can be regarded as a core of a 2-handle for $C_2$. Then by taking
the handle decomposition for $M$ obtained by $\mathcal{D}_1$ and
$\mathcal{D}_2$ we obtain the conclusion of Lemma~\ref{lem:haken
annulus disjoint from handles}.
\end{proof}

Suppose that $M$ is the union of two compact orientable manifolds
(say $M = M^{(1)} \cup M^{(2)}$) such that $M^{(1)} \cap M^{(2)}=
\del M^{(1)} \cap \del M^{(2)} = F$, a collection of components of
$\del M^{(i)}$.  Let $C_1^{(1)} \cup_{\s^{(1)}} C_2^{(1)}$ and
$C_1^{(2)} \cup_{\s^{(2)}} C_2^{(2)}$ be Heegaard splittings of
$M^{(1)}$ and $M^{(2)}$ respectively.  Let $F^{(1)} \cup F_1^{(1)}
\cup F_2^{(1)}$ and $F^{(2)} \cup F_1^{(2)} \cup F_2^{(2)}$ be
partitions of the components of $\del M^{(1)}$, $\del M^{(2)}$
respectively such that the surfaces $F^{(1)}$ and $F^{(2)}$ are
identified in $M$ as $F$, and

$$\del_- C_1^{(1)} = F_1^{(1)},\mbox{\hspace{.5in}}
\del_- C_2^{(1)} = F^{(1)} \cup F_2^{(1)},$$

$$\del_- C_1^{(2)} = F^{(2)} \cup F_1^{(2)},\mbox{\hspace{.5in}}
\del_- C_2^{(2)} = F_2^{(2)}.$$

Let $C_1 \cup_\s C_2$ be an amalgamation of  $C_1^{(1)}
\cup_{\s^{(1)}} C_2^{(1)}$ and $C_1^{(2)} \cup_{\s^{(2)}}
C_2^{(2)}$.

\begin{pro}
\label{pro:haken annulus after amalgamation}

Let $M$, $M^{(1)}$, $M^{(2)}$,  $C_1 \cup_\s C_2$, $C_1^{(1)}
\cup_{\s^{(1)}} C_2^{(1)}$ and $C_1^{(2)} \cup_{\s^{(2)}}
C_2^{(2)}$ be as above.  Suppose that there exists a Haken annulus
$A^{(1)}$ for $C_1^{(1)} \cup_{\s^{(1)}} C_2^{(1)}$ such that
$A^{(1)} \cap F^{(1)} = \emptyset$. Then the image of $A^{(1)}$ is a
Haken annulus for $C_1 \cup_\s C_2$.
\end{pro}

\begin{proof}

By Lemma~\ref{lem:haken annulus disjoint from handles} there
exists a handle decomposition of $M^{(1)}$ of the form

$$ M^{(1)} = (\del_- C_1^{(1)} \times [0,1])
\cup \mathcal{H}_1^{(1)}
  \cup \mathcal{H}_2^{(1)}$$

\noindent induced by $C_1^{(1)} \cup_{\s^{(1)}} C_2^{(1)}$ such
that $A^{(1)} \cap (\mathcal{H}_1^{(1)} \cup \mathcal{H}_2^{(1)})
= \emptyset$. By the definition of amalgamation, we can take a
handle decomposition:

$$M^{(2)} = (F^{(2)} \times [0,1]) \cup (F_1^{(2)} \times [0,1])
 \cup \mathcal{H}_0^{(2)}
\cup \mathcal{H}_1^{(2)}
  \cup \mathcal{H}_2^{(2)} \cup \mathcal{H}_3^{(2)}$$

\noindent induced by $C_1^{(2)} \cup_{\s^{(2)}} C_2^{(2)}$ such
that these handles give the handle decomposition:

$$M = (F_1^{(1)} \times [0,1]) \cup (F_1^{(2)} \times [0,1])
 \cup \mathcal{H}_0^{(2)}
\cup (\mathcal{H}_1^{(1)}\cup \mathcal{H}_1^{(2)})\cup$$

$$\cup (\mathcal{H}_2^{(1)} \cup \mathcal{H}_2^{(2)})
  \cup \mathcal{H}_3^{(2)},$$

\noindent which gives the Heegaard splitting $C_1 \cup_\s C_2$.
Here we note that adding $F_1^{(2)} \times [0,1]$ and
$\mathcal{H}_i^{(2)}$ does not affect $A^{(1)}$, since $A^{(1)}
\cap F = \emptyset$.  Hence we see that $A^{(1)} \cap (\cup_{i,j}
\mathcal{H}_i^{(j)}) = \emptyset$, and this together with
Lemma~\ref{lem:obvious Haken annulus} shows that the image of
$A^{(1)}$ is a Haken annulus for $C_1 \cup_{\s} C_2$.
\end{proof}

\begin{pro}
\label{pro:combined haken annuli} Let $M$, $M^{(1)}$, $M^{(2)}$,
$C_1 \cup_\s C_2$, $C_1^{(1)} \cup_{\s^{(1)}} C_2^{(1)}$, and
$C_1^{(2)} \cup_{\s^{(2)}} C_2^{(2)}$ be as above.  Suppose that
there exists a Haken annulus $A^{(j)}$ for $C_1^{(j)}
\cup_{\s^{(j)}} C_2^{(j)}$ for each $j=1,2$, such that $A^{(1)}
\cap F = A^{(2)} \cap F$ (hence the image of $A^{(1)} \cup
A^{(2)}$ is an annulus properly embedded in $M$).
Then the image of $A^{(1)} \cup A^{(2)}$ is a Haken annulus for
$C_1 \cup_\s C_2$.
\end{pro}

\begin{proof}

By Lemma~\ref{lem:haken annulus disjoint from handles} there
exists a handle decomposition $(F_1^{(1)} \times [0,1]) \cup
\mathcal{H}_1^{(1)}\cup \mathcal{H}_2^{(1)}$ ($(F_1^{(2)} \times
[0,1]) \cup (F^{(2)} \times [0,1]) \cup \mathcal{H}_1^{(2)}\cup
\mathcal{H}_2^{(2)}$ resp.) induced by $C_1^{(1)} \cup_{\s^{(1)}}
C_2^{(1)}$ ($C_1^{(2)} \cup_{\s^{(2)}} C_2^{(2)}$ resp.) such that
$A^{(1)} \cap (\mathcal{H}_1^{(1)} \cup \mathcal{H}_2^{(1)}) =
\emptyset$ ($A^{(2)} \cap (\mathcal{H}_1^{(2)} \cup
\mathcal{H}_2^{(2)}) = \emptyset$ resp.). By pushing $F^{(2)}
\times [0,1]$ into $M^{(1)}$, we obtain a handle decomposition $M=
(F_1^{(1)} \times [0,1]) \cup (\mathcal{H}_1^{(1)} \cup
\mathcal{H}_2^{(1)}) \cup (F_1^{(2)} \times [0,1]) \cup
(\mathcal{H}_1^{(2)} \cup \mathcal{H}_2^{(2)})$ (see
Section~\ref{sec:background}). Here we may suppose that the image
of $A^{(1)} \cup A^{(2)}$ is 
$A^{(1)}$. 
It is easy to see that by isotopy we may suppose,
in addition to the above, that $\mathcal{H}_1^{(2)} \cap
(\mathcal{H}_1^{(1)} \cup \mathcal{H}_2^{(1)})=\emptyset$ and
$\mathcal{H}_2^{(2)} \cap \mathcal{H}_2^{(1)} = \emptyset$.  Hence
these handles give the Heegaard splitting $C_1 \cup_\s C_2$.  For
these handles we have:

$$A^{(1)} \cap
(\mathcal{H}_1^{(1)} \cup \mathcal{H}_1^{(2)} \cup
\mathcal{H}_2^{(1)} \cup \mathcal{H}_2^{(2)})= \emptyset.$$

Then Lemma~\ref{lem:obvious Haken annulus} establishes the
conclusion of Proposition~\ref{pro:combined haken annuli}.
\end{proof}

\section{The swallow follow torus theorem}
\label{sec:sft}

Let $X$ be the exterior of a knot $K$ in a closed orientable
3-manifold. Let
$\gamma_1,\dots,\gamma_p$ be mutually disjoint \scc s that are
obtained by pushing mutually disjoint meridional curves in $\del
X$ into $\mbox{\rm int}(X)$.   Then let $X^{(p)} = \mbox{\rm cl}(X
\setminus \cup_{i=1}^p N(\gamma_i))$.  Set $X^{(0)} = X$. A torus $T$ in
$X^{(p)}$ is called a \em\sft\ \em if it is separating and essential in
$X^{(p)}$, and there exists an annulus $A$ in $X^{(p)}$ so that $A
\cap \del X^{(p)} = \del A \cap \del X^{(p)}$ is a meridian curve
in $\del X$ and $A \cap T = \del A \cap T$ is an essential \scc\
in $T$.  We call this annulus a \em meridionally compressing
annulus \em for $T$.  Let $A'$ be the annulus properly embedded in $X^{(p)}$
obtained from $T$ by removing $N(A) \cap T,$ then adding two parallel copies
of $A$.  (We say that $A'$ is obtained from $T$ by performing surgery along
$A$.)  Since $T$ is essential we see that $A'$ is also essential.
The annulus $A'$ gives a decomposition:

$$X_1^{(p_1)} \cup_{A'} X_2^{(p_2)}, \mbox{\hspace{.2in}} p=p_1+p_2,$$

\noindent where each $X_i$ is the exterior of a knot $K_i$
(possibly the trivial knot in $S^3$, in which case $X_i^{(p_i)}$
is homeomorphic to (disk with $p_i$ holes) $\times S^1$; recall
that this was denoted $D(p_i)$ at the end of
Section~\ref{sec:background}). We call this decomposition the \em
decomposition induced by $T$ \em. Note that $T$ can be retrieved
from $A'$ by adding a component of $\del X \setminus \del A'$ and
applying an isotopy pushing it into $\mbox{\rm int} X^{(p)}$. The
component of $\del X \setminus \del A'$ follows $K_1$ or $K_2$,
say $K_2$. This yields the following decomposition:

$$X^{(p)} = X_1^{(p_1 + 1)} \cup_T X_2^{(p_2)}.$$

The purpose of this section is to prove the following:

\begin{thm}[The Swallow Follow Torus Theorem]
\label{thm:sft}

Let $K_1,\dots,K_n$ ($n\geq 2$) be knots in closed orientable
3-manifolds.  Let $X_1,\dots,X_n$ and $X$ be the exteriors of
$K_1, \dots, K_n$ and $\#_{i=1}^n K_i$ respectively. Suppose that
no $X_i$ is a solid torus and $X$ is irreducible.  Then:

\begin{enumerate}
    \item For any $p \geq 0$, every Heegaard surface for
    $X^{(p)}$ weakly reduces to exactly one surface which is a \sft, or---
    \item For some $i$, $X_i$ contains an essential meridional surface; more
    precisely, if a genus $g$ \hhs\ for $X^{(p)}$ (for some $p$) does not
    weakly reduce to a \sft, then
    there exists a meridional essential surface $S \subset X_i$ so
    that $\chi(S) \geq 4-2g(X^{(p)})$.
\end{enumerate}
\end{thm}

\begin{rmkk}{\rm
\label{rmk:no sft when n 1} One might hope that Theorem
\ref{thm:sft} holds for $n=1$.  However, this is not the case. Let
$K$ be a non-trivial 2-bridge knot in $S^3$.  Then $K$ is m-small
(\cite{hatcher-thurston}) and $g(X) = g(X^{(1)})=2$ 
(see, for example, \cite{morimoto1994}). 
If a genus 2 Heegaard surface
weakly reduces then the obtained essential surface is a 2-sphere.
Applying this to $X^{(1)}$, we see that (since $X^{(1)}$ is
irreducible) every minimal genus Heegaard splitting of $X^{(1)}$
is strongly irreducible; hence the conclusion of
Theorem~\ref{thm:sft} does not hold for $X^{(1)}$. }\end{rmkk}

As an immediate consequence of Theorem~\ref{thm:sft} we have:

\begin{cor}
\label{cor:sft}

Let $K_i$, $K$, $X_i$ and $X$ be as in the statement of Theorem
\ref{thm:sft}, and suppose that each $K_i$ is m-small.  Then for
any $p \geq 0,$ every Heegaard splitting of $X^{(p)}$ weakly
reduces to a \sft.
\end{cor}

\begin{proof}[Proof of Theorem~\ref{thm:sft}]

Note that if some $K_i$ is a composite knot then $X_i$ contains an
essential meridional annulus and conclusion (2) of
Theorem~\ref{thm:hopf-annulus} holds. We assume from now on that
for each $i$, $K_i$ is prime.

\begin{clm}
\label{clm:number of factors} Either existence and uniqueness of
prime factorization holds for $K$ or conclusion (2) holds.
\end{clm}

We note that prime factorization of knots in 3-manifolds is
studied by K. Miyazaki in \cite{miyazaki} and it is shown that
neither existence nor uniqueness holds in general.

\begin{proof}[Proof of Claim~\ref{clm:number of factors}]
It is shown in the Existence Theorem of \cite{miyazaki} that if
the meridian of a knot $K$ is non-null homotopic then the number
of prime factors of $K$ is unique.  If the meridian of $K$ is null
homotopic then by the Loop Theorem (\eg\ \cite{hempel}) $\del X$
is compressible, and either $X$ is a solid torus or it is
reducible, both contrary to our assumptions. The Uniqueness
Theorem of \cite{miyazaki} implies that if the uniqueness of prime
decomposition of $K$ fails then the exterior of some factor of $K$
contains a non-separating meridional annulus; this gives
conclusion (2).
\end{proof}

By Claim~\ref{clm:number of factors}, for the remainder of the
proof we may suppose that $K$ has a unique prime factorization.

Note that if a given Heegaard splitting weakly reduces to a \sft\
then every stabilization of that Heegaard splitting can be weakly reduced
to the same \sft; it therefore suffices to prove the theorem for
non-stabilized Heegaard splittings.  Fix $p$ and let $C_1 \cup_\s
C_2$ be a non-stabilized \hh\ splitting of $X^{(p)}$.

Let $\mathcal{A} = A_1 \cup \cdots \cup A_n$ be a union of mutually disjoint
annuli properly embedded in $X^{(p)}$ with the following properties:
\begin{enumerate}
\item For each $i$, $\del A_i \subset \del X$,
\item $X^{(p)}$ cut open along $\mathcal{A}$ consists of $n+1$ components, $n$
  homeomorphic to $X_1,\dots,X_n$ (by slightly abusing notation $X_i$ will
  refer to the component homeomorphic to $X_i$) and one component homeomorphic
  to $D(p)$ (as defined in the end of Section~\ref{sec:background}), denoted $X_0$,
  and:
\item For each $i$, $X_i \cap \mathcal{A} = A_i$.
\end{enumerate}
We say that the annuli of $\mathcal{A}$ are {\it not nested}. Let
$T_j$ ($j=1,\dots,p$) denote the component of $\del X^{(p)}$
corresponding to $\del N(\gamma_j).$ (Hence $\del X^{(p)} = \del X
\cup (\cup_{j=1}^p T_j)$.)

We divide the proof into the following cases.

{\bf Case 1. $\s$ is strongly irreducible.} Since $\s$ is strongly
irreducible, we can isotope it to intersect $\mathcal{A}$ in \scc
s that are essential in both $\mathcal{A}$ and $\s$.  Minimize
$|\mathcal{A} \cap \s|$ subject to that constraint. Let $\s_i = \s
\cap X_i$ ($i=0,1,\dots,n$). Suppose that $\s_i = \emptyset$ for
some $i$ ($1 \le i \le n$). Then $A_i$ is contained in a
compression body $C_j$ ($j=1$, or $2$). By (6) of
Remarks~\ref{rmk:properties of CBs}, we see that $A_i$ is boundary
parallel in $C_j$. This shows that $X_i$ is a solid torus,
contradicting the assumption of Theorem~\ref{thm:sft}. Hence $\s_i
\ne \emptyset$ for each $i$ ($i=1, \dots , n$).

The proof of the following claim is similar to the proof of Claim~4 in page
247 of \cite{kobayashi-rieck-local-det}, however for the convenience of the
reader we include it here:

\begin{clm}
\label{clm:s k incomprresible}
\begin{enumerate}
\item There exists $i$ (possibly $i=0$) such that $\Sigma_i$ is compressible
  in $X_i$, and there exist two compressing disks $D_1$, $D_2 \subset X_i$
  that compress $\s_i$ into the $C_1$ side and $C_2$ side respectively, and
  $D_1 \cap D_2 \neq \emptyset$.
\item For each $k \neq i$, $\s_k$ is incompressible in $X_k$.
\end{enumerate}
\end{clm}

\begin{proof}[Proof of Claim 2]
Note that since $\s \cap \mathcal{A}$ is essential in $\mathcal{A}$ any
compressing disk for $\s_i$ is also a compressing disk for $\s$.

Subclaim: For each $i=1,2$ there is a meridional disk $D_i$ of $C_i$ so that
$D_i \cap \mathcal{A} = \emptyset$.

\noindent Proof.   Fix $i$ and let $D_i \subset C_i$ be a meridian
disk that minimizes $|D_i \cap \mathcal{A}|$.  If $D_i \cap
\mathcal{A} = \emptyset$ we are done.  Else, by minimality of
$|D_i \cap \mathcal{A}|$ and a standard innermost disk argument,
no component of $D_i \cap \mathcal{A}$ is a \scc, and each
outermost disk in $D_i$ gives a boundary compression of $\s_k$ in
$X_k$ (for some $k$).  By Lemma~\ref{boudary compressible implies
compressible} we see that either there is a compression of $\s_k$
in $X_k$ which is contained in $C_i$, or a component of $\s_k$ is
a boundary parallel annulus in $X_k$; the former is the desired
conclusion of the subclaim, and the latter contradicts minimality
of $|D_i \cap \mathcal{A}|$, establishing the subclaim.

Let $D_1$, $D_2$ be as in Subclaim.  If $D_1 \cap D_2 =
\emptyset$ then $\{D_1,D_2\}$ gives a weak reduction of $\s$,
contrary to our assumption. Hence $D_1$ and $D_2$ are contained in
the same component (say $X_i$). Finally note that if for $k \neq
i$ $\s_k$ compresses in $X_k$ then the compressing disk for $\s_k$
and either $D_1$ or $D_2$ gives a weak reduction of $\s$, a
contradiction; hence for each $k \neq i$, $\s_k$ is
incompressible.
\end{proof}

Since $n \geq 2,$ we may suppose that $\s_1$ satisfies (2), i.e.,
$\s_1$ is incompressible in $X_1$.  Let $\s'$ be the surface
obtained by compressing $\s$ along $D_1$.   Since $\mathcal{A}$ is
incompressible and $\mathcal{A} \cap \s$ is essential in
$\mathcal{A}$, no component of $\s'$ cut open along $\mathcal{A}$
has positive Euler characteristic.  Hence $\chi(\s_1) = \chi(\s'
\cap X_1) \geq 4-2g$. Suppose that each component of $\s_1$ is
boundary parallel in $X_1$. Then the minimality of $|\mathcal{A}
\cap \s|$ implies that each component of $\s_1$ is parallel in
$X_1$ to the annulus $X_1 \cap \del X$. Let $X_1' = \mbox{cl}(X_1
\setminus N(\del X_1,X_1))$. Then we can push $\s \cap X_1'$ out
of $X_1'$, hence, we may assume $X_1' \cap \s = \emptyset$, and
$X_1'$ is contained in a compression body $C_j$ ($j=1$, or $2$).
Note that $\del X_1'$ is incompressible in $X$, and therefore in
$C_j$.  Hence $\del X_1'$ is parallel to $\del_- C_j$ ((6) of
Remarks~\ref{rmk:properties of CBs}).   This shows that $X^{(p)}$
is homeomorphic to $X_1'$, but this contradicts the assumption
that $n \geq 2$.  Hence some component of $\s_1$ is not boundary
parallel in $X_1$. This establishes conclusion~(2) of
Theorem~\ref{thm:sft}.

{\bf Case 2. $\s$ is weakly reducible.} Take a
Scharlemann--Thompson untelescoping (Section~\ref{sec:background})
of $C_1 \cup_\s C_2$ to obtain a (possibly disconnected)
incompressible surface $S$ (corresponding to the union of
$\partial R_j$ in Remark~\ref{rmk:scharlemann's ST-untel is
different} minus $\partial X^{(p)}$) with $\chi(S) \geq 6-2g$. 
Denote the connected
components of $S$ by $S_j$. Since $C_1 \cup_\s C_2$ is
non-stabilized and $X$ is irreducible, we see, by
Proposition~\ref{pro:untel to a conn sep surfce implies weak
reduction}, Lemma~\ref{lem:weak reduction that yeilds a 2-sphere},
and Waldhausen's classification of Heegaard surfaces of $S^3$
\cite{waldhausen}, that no $S_j$ is a 2-sphere. Hence, for each
$j$, $\chi(S_j) \geq 6-2g$. Since both $S$ and $\mathcal{A}$ are
incompressible, we may suppose that each component of $S \cap
\mathcal{A}$ is a \scc\ which is essential in both $S$ and
$\mathcal{A}$.  Minimize $|S \cap \mathcal{A}|$ subject to this
constraint.

We have the following subcases.

\noindent {\bf Subcase 2a. $S \cap \mathcal{A} = \emptyset.$} Let
$M$ be the the component of $X^{(p)}$ cut open along
$S$ containing $\del X$. Note that in this case $\mathcal{A}
\subset M$. Let $C_{M,1} \cup_{\s_M} C_{M,2}$ be the Heegaard
splitting of $M$ induced from $\s$.  Since $\mathcal{A} \subset
M$, $M$ is not a product $\del X \times [0,1]$.  Since $C_{M,1}
\cup_{\s_M} C_{M,2}$ is strongly irreducible, we may suppose each
component of $\s_M \cap \mathcal{A}$ is essential in both $\s_M$
and $\mathcal{A}$; minimize $|\s_M \cap \mathcal{A}|$ subject to
this constraint.

Then we claim that $\s_M \cap \mathcal{A} \neq \emptyset$.
Suppose, for a contradiction, that $\s_M \cap \mathcal{A} = \emptyset$.
Then by (6) of Remarks~\ref{rmk:properties of CBs} we have that
each $A_i$ is boundary parallel in $C_{M,1}$.  Hence each $A_i$ is
boundary parallel in $X^{(p)}$.  This shows that each $X_i$ is a
solid torus, contradicting the assumption of Theorem~\ref{thm:sft}.

By Claim~\ref{clm:s k incomprresible} above we may assume that
each component of $\s_M \cap (X_1 \cap M)$ (=$\s_M \cap X_1$) is
incompressible in $X_1 \cap M$. Note that each component of
$\mbox{\rm Fr}_{X^{(p)}}(X_1 \cap M)$ is either $A_1$ or a
component of $S$, and therefore is incompressible in $X_1$. This
shows that each component of $\s_M \cap X_1$ is incompressible in
$X_1$.

If some component of $\s_M \cap X_1$ is not boundary parallel in
$X_1$, then we have conclusion (2) of Theorem~\ref{thm:sft}.
Suppose that each component of $\s_M \cap X_1$ is a boundary
parallel annulus in $X_1$.  Pick a component of $\s_M \cap X_1$
and let $V$ be the product region that it separates from $X_1$; note
that $V$ is a solid torus. Since each component of $S$ is
incompressible in $X$, we see that no component of $S$ is
contained in $V$.  
Let $X_1' = \mbox{cl}(X_1 \setminus N(\del X_1,X_1))$. 
Since $S \cap V = \emptyset$ we can push $\s_M
\cap X_1'$ out of $X_1'$. Hence we have $\del X_1' \subset
C_{M,1}$ or $\del X_1' \subset C_{M,2}$, say $C_{M,1}$. 
Since 
$\del X_1'$ is incompressible in $X$, it is incompressible in
$C_{M,1}$. Hence $\del X_1'$ is parallel to a component of $\del_-
C_{M,1}$ ((6) of Remarks~\ref{rmk:properties of CBs}). Note that
each component of $\del_- C_{M,1}$ is either $\del X$ or a
component of $S$.  Since $X \neq X_1$ we see that $\del X_1'$ is
parallel to a component of $S$.  By Proposition~\ref{pro:untel to
a conn sep surfce implies weak reduction} there is a weak
reduction of $C_1 \cup_\s C_2$ which gives exactly $\del X_1'$.
This establishes conclusion~(1) of Theorem~\ref{thm:sft}.

{\bf Subcase 2b.  $S \cap \mathcal{A} \neq \emptyset$.}

Suppose that the conclusion (2) of Theorem~\ref{thm:sft} does not
hold, hence, for each $i \neq 0$, each component of $S \cap X_i$ is
closed, or a boundary parallel annulus.

Let $X_i' = \mbox{cl}(X_i \setminus N(\del X_i,X_i))$
($i=1,\dots,n$) and $X_0' =\mbox{cl}(X \setminus \cup_{i=1}^n
X_i')$. Note that $X_0'$ is homeomorphic to a (disk with $n+p$
holes) $\times S^1$.  Since each component of $S \cap X_i$ with
non-empty boundary is a boundary parallel annulus, they can be
pushed into $X_0'$. Hence we have $S \cap \del X_0' = \emptyset$,
and $S \cap X_0' \ne \emptyset$. Since $S$ is a closed
incompressible surface in the Seifert fibered space $X_0'$ it must
be vertical by \cite[Theorem VI.34]{jaco}, that is to say, $S$ is
isotopic to a surface of the form $\mbox{(simple closed curve(s))}
\times S^1$. Note that any essential \scc\ on a disk with $n+p$
holes is separating. Hence each component of $S \cap X_0'$ is
separating in $X_0'$, and therefore also in $X^{(p)}$. Suppose that
there is a component of $S \cap X_0'$, say $S_0$, which is not
boundary parallel in $X^{(p)}$. Then $S_0$ separates $X_0'$ into
two components, and let $X_0''$ denote the closure of the
component of $X_0 \setminus S_0$ such that $\del X \subset X_0''$.
Let $A$ be a vertical annulus properly embedded in $X_0''$ such
that one component of $\del A$ is contained in $\del X$, and the
other in $S_0$. Note that $A \cap \del X$ is a meridian curve and
$A \cap S_0$ is an essential \scc\ in $S_0$. Hence $S_0$ is a
\sft \ in $X^{(p)}$.  
Finally, by Proposition~\ref{pro:untel to a conn sep surfce
implies weak reduction}, we see that a weak reduction of $C_1
\cup_P C_2$ gives $S_0$. This establishes conclusion~(1) of
Theorem~\ref{thm:sft}.

Finally suppose that each component of $S \cap X_0'$ is boundary
parallel in $X^{(p)}$. Let $\widetilde{P}$ be the union of the
parallel regions bounded by the components of $S \cap X_0'$. Then
let $\widetilde{X}^{(p)} = \mbox{cl}(X^{(p)} \setminus
\widetilde{P})$, and $\widetilde{C}_1
\cup_{\tilde{\Sigma}}\widetilde{C}_2$ be the Heegaard splitting of
$\widetilde{X}^{(p)}$ obtained from the Scharlemann-Thompson
untelescoping of $C_1 \cup_\Sigma C_2$ by amalgamating the pieces
contained in $\widetilde{X}^{(p)}$. Then we apply the above
argument to $\widetilde{C}_1 \cup_{\widetilde{\s}}
\widetilde{C}_2$, that is, we apply the argument of Case 1 if
$\widetilde{\s}$ is strongly irreducible or the argument of Case 2
is $\widetilde{\s}$ is weakly reducible.  Note that for
$\widetilde{\s}$ we do not have the case corresponding to Subcase
2b with each component of $S \cap X_0'$ being boundary parallel in
$X^{(p)}$.  Then we have either conclusion~(1) or conclusion~(2)
for $\widetilde{C}_1 \cup_{\widetilde{\s}} \widetilde{C}_2$.
Conclusion~(2) for $\widetilde{\s}$ immediately establishes
conclusion~(2) for $C_1 \cup_\Sigma C_2$. If we have
conclusion~(1) for $\widetilde{C}_1
\cup_{\tilde{\Sigma}}\widetilde{C}_2$, then this together with
Proposition~\ref{pro:untel to a conn sep surfce implies weak
reduction} establishes conclusion~(1) for $C_1 \cup_\Sigma C_2$.

This completes the proof of Theorem~\ref{thm:sft}.
\end{proof}

\section{Numerical bounds}
\label{sec:nimerical-bounds}

\begin{thm}
\label{thm:numerical-bounds}

Let $X$ be the exterior of a knot $K$ in a closed orientable
3-manifold. Suppose that there exists a minimal genus Heegaard
surface for $X$ that weakly reduces to a swallow follow torus
inducing the decomposition $K = K_1 \# K_2$; let $X_i = E(K_i)$.
Then we have:

$$g(X_1) + g(X_2) -1 \leq g(X) \leq g(X_1) + g(X_2).$$
\end{thm}

\begin{rmkk}
{\rm Recall from Section~\ref{sec:introduction} that the right
hand side inequality always holds.}
\end{rmkk}

As a consequence of Theorem~\ref{thm:numerical-bounds} and
Corollary~\ref{cor:sft}, we have:

\begin{cor}
\label{cor:bounding-degeneration}

Let $K_1,\dots,K_n$ be m-small knots in closed orientable
3-manifolds, let $X_i = E(K_i)$, and $X = E(\#_{i=1}^n K_i)$.
Suppose that no $X_i$ is a solid torus and $X$ is irreducible.
Then we have the following inequalities:

$$\Sigma_{i=1}^n g(X_i) - (n-1) \leq g(X) \leq \Sigma_{i=1}^n
g(X_i).$$
\end{cor}

\begin{proof}[Proof of Theorem~\ref{thm:numerical-bounds}]

By Section~\ref{sec:sft} we have a decomposition: $X = X_1^{(1)}
\cup X_2$ or $X=X_1 \cup  X_2^{(1)}$.  Since the argument is
symmetric we may suppose that $X=X_1^{(1)} \cup X_2$.

Note that $X_1$ is obtained from $X_1^{(1)}$ by Dehn filling and
that the core of the attached solid torus is parallel to a
meridian curve on $\del X_1$.
Since any Heegaard surface for $X_1$ can be regarded as the
frontier of a compression body obtained from $N(\del X,X)$ by
attaching 1-handles, this implies that the core of the attached
solid torus, say $l$, is isotopic into every Heegaard surface of
$X_1$, that is to say, the filling is called {\it good}, as in
\cite[Definition 4.1]{rieck} (see also
\cite{rieck-sedgwick2}\cite{rieck-sedgwick1}).  Let $C_1'
\cup_{\s'} C_2'$ be a minimal genus Heegaard splitting of $X_1$
such that $\del X_1 \subset C_1'$, and $l \subset C_2'$. Let
$\widetilde{\s}'$ be a surface obtained from $\s$ and $\del N(l)$ 
as is shown in \cite[Corollary 4.3]{rieck}; 
note that $g(\widetilde{\s}') =
g(\s') + 1$ (in fact, $\widetilde{\s}'$ is a stabilization of
$\s'$) and $\widetilde{\s}'$ is a Heegaard surface for
$X_1^{(1)}$.  Thus we get the inequality $g(X_1^{(1)}) \leq g(X_1)
+ 1$. Since $X_1$ is obtained from $X_1^{(1)}$ by Dehn filling, we
immediately see that $g(X_1) \leq g(X_1^{(1)})$. Combining the two
we have:

$$g(X_1) \leq g(X_1^{(1)}) \leq g(X_1) +1.$$

\noindent
On the other hand, by Proposition~\ref{prop:amalgamation of
minimal genus: connected case} we have the following equality:

$$g(X) = g(X_1^{(1)}) + g(X_2) - 1.$$

\noindent
Therefore we have:

$$g(X_1) + g(X_2) -1 \leq g(X) \leq g(X_1) + g(X_2).$$
\end{proof}

In Section~\ref{sec:morimoto}, we will need the following information
which follows from the proof of Theorem~\ref{thm:numerical-bounds},
and so we single it out here.

\begin{cor}
\label{cor:when can the genus go down}

If a minimal genus Heegaard splitting for $X$ weakly reduces to a
swallow follow torus $T$ inducing $X = X_1^{(1)} \cup_T X_2$, then
either $g(X_1) = g(X_1^{(1)})$ or $g(X_1) = g(X_1^{(1)}) - 1$, and
exactly one of the following holds:

\begin{enumerate}
    \item $g(X) = g(X_1) + g(X_2)$ and $g(X_1) = g(X_1^{(1)}) - 1$.
    \item $g(X) = g(X_1) + g(X_2) -1$ and $g(X_1) = g(X_1^{(1)})$.
\end{enumerate}
\end{cor}

\section{Hopf-Haken Annulus Theorem}
\label{sec:hopf-annulus}

Let $K_1,\dots,K_n$ ($n\geq 1$) be knots in closed orientable
3-manifolds.  Let $X_1,\dots,X_n$, and $X$ be the exteriors of
$K_1, \dots, K_n$, and $K = \#_{i=1}^n K_i$ respectively.

Recall from Section~\ref{sec:sft} that for $p \geq 1$, $X^{(p)} =
\mbox{\rm cl}(X \setminus \cup_{i=1}^p N(\gamma_i))$, where
$\gamma_1, \dots , \gamma_p$ are \scc s obtained by pushing
mutually disjoint meridian curves in $\del X$ into $\mbox{int}X$
(with $X^{(0)} = X$).  One can think of $X^{(p)}$ as the exterior
of the link obtained from $K$ by taking connected sum with the
Hopf link $n$ times. Note that the exterior of the Hopf link has
an essential annulus connecting the distinct boundary components,
where one boundary component of the annulus is a meridian curve,
and the other is a longitude. This shows that $X^{(p)}$ has an
essential annulus as well (provided $p \geq 1$) and motivates the
following definition (\cf\ Definition~\ref{dfn:haken annulus}):

\begin{dfn}[Hopf-Haken annulus]
\label{dfn:hopf haken annulus}

A \em Hopf spanning annulus \em in $X^{(p)}$ is a properly
embedded essential annulus $A$ with one component of $\del A$ a
meridian curve of $\del X$ and the other component a longitude of
$\del N(\gamma_i)$ for some $i$. If in addition there exists a
Heegaard splitting $C_1 \cup_{\s} C_2$ of $X^{(p)}$ such that $A$
intersects $\s$ in one simple closed curve that is essential in
$A$, then $A$  is called a \em Hopf-Haken annulus \em for $C_1
\cup_{\s} C_2$. In that case, $\s$ is said to admit a Hopf-Haken
annulus.
\end{dfn}

By Figure~6.1, we immediately have the following.

\begin{figure}[ht]
\begin{center}
\includegraphics[width=9cm, clip]{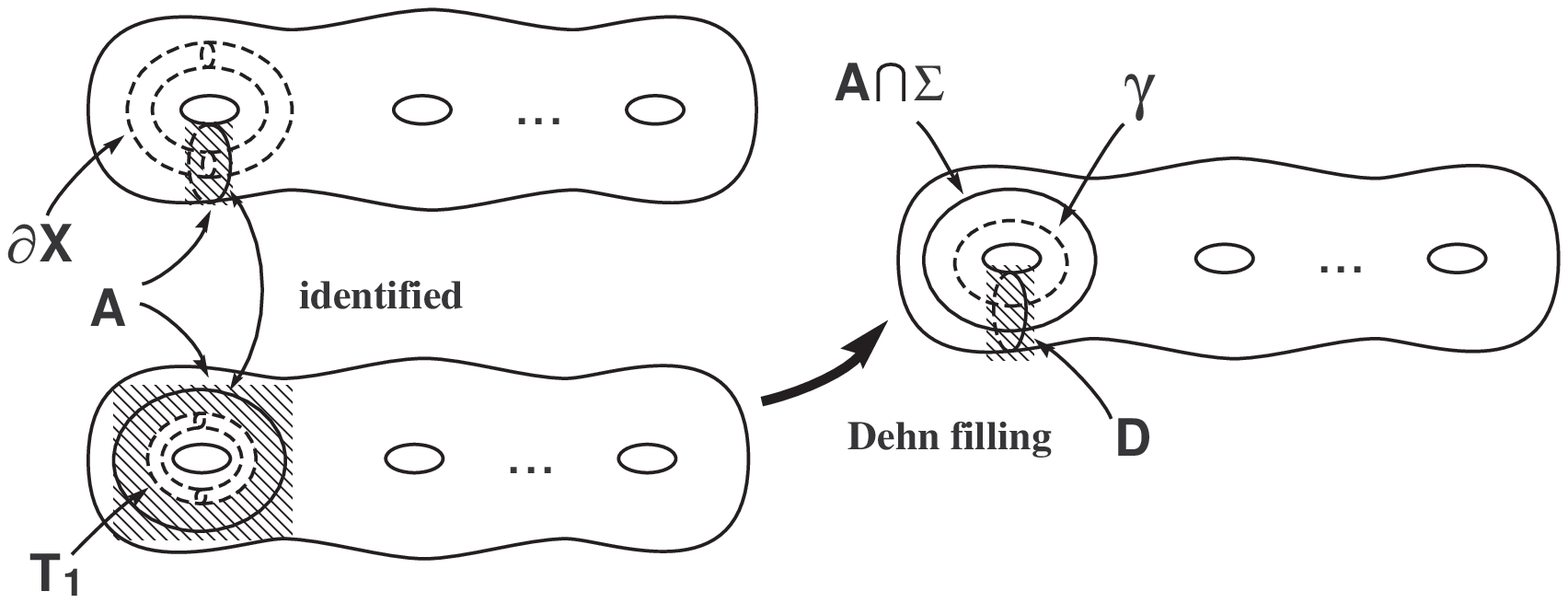}
\end{center}
\begin{center}
Figure 6.1.
\end{center}
\end{figure}

\begin{pro}
\label{pro:hopf-haken implies primitive meridian}

Let $Y$ be the exterior of a knot $K$.  Suppose that for some $p
\geq 1$ there exists a Heegaard surface $\s$ of
$Y^{(p)}$ such that $g(\s) = g(Y)$ and $\s$ admits a Hopf--Haken
annulus.  Then $K$ admits a primitive meridian.
\end{pro}

In this section we prove the following theorem.

\begin{thm}[The Hopf-Haken Annulus Theorem]
\label{thm:hopf-annulus} Let $K_i$ ($i=1,\dots,n$, $n \geq 1$) be
knots in closed orientable 3-manifolds.  Let $X_i$, and $X$ be the
exteriors of $K_i$, and $K = \#_{i=1}^n K_i$ respectively. Let
$X^{(p)}$ be as above. Suppose that no $X_i$ is a solid torus and
$X$ is irreducible. Then one of the following holds.
\begin{enumerate}
    \item For each $p$ $(\ge 1)$,
    there exists a minimal genus Heegaard surface of
    $X^{(p)}$ that admits a Hopf-Haken annulus.
    \item For some $i$ ($1 \leq i \leq n$),
    $X_i$ contains a meridional essential surface $S$ such that
    $\chi(S) \geq 4-2g$,
    where $g = g(X^{(p)})$.
\end{enumerate}
\end{thm}

\begin{rmk}

Note that in Theorem~\ref{thm:hopf-annulus} we cannot expect every
minimal genus Heegaard surface to satisfy conclusion~(1), even
when conclusion~(2) is not satisfied.  To see this, consider the
following example: let $K_m$ be the knots constructed by
Morimoto--Sakuma--Yokota in \cite{skuma-morimoto-yokata} and $X_m$
their exteriors. Then $g(X_m)=2$ (by
\cite{skuma-morimoto-yokata}) and $g(X_m^{(1)})=3$ (by
\cite{kobayashi-rieck-growth-rate}). There is an essential torus
$T \subset X_m$ giving the decomposition $X_m^{(1)} = D(2) \cup
X_m'$, where $X_m'$ is homeomorphic to $X_m$, 
and $D(2)$ is (disk with two holes)$\times S^1$. 
Note that $\partial
D(2) =
\partial X_m^{(1)} \cup T$. By Proposition~\ref{pro:genus of D(n)}
and the comment after it, $g(D(2))=2$ and there exists a genus two
Heegaard splitting of $D(2)$ (say $C_1' \cup_{\s'} C_2'$) with
$\partial_-C_1' = \partial X_m^{(1)}$ and $\partial_-C_2' = T$.
Let $C_1 \cup_\s C_2$ be the Heegaard splitting 
of $X_m^{(1)}$ obtained by
amalgamating $C_1' \cup_{\s'} C_2'$ and a minimal genus Heegaard
splitting of $X_m'$. By Lemma~\ref{lem:genus after amalgamation},
$C_1 \cup_\s C_2$ is a minimal genus Heegaard splitting of
$X_m^{(1)}$. By the definition of amalgamation we see that
$\partial_- C_1 = \partial X_m^{(1)}$.  We conclude that $C_1
\cup_\s C_2$ does not admit a Hopf--Haken annulus, for if it did
then the components of $\partial X_m^{(1)}$ would be in distinct
compression bodies.

Although we do not know if $K_m$ is m-small, this example seems to
indicate that replacing ``there exists" with ``for every" in
conclusion~(1) of Theorem~\ref{thm:hopf-annulus} is not likely to
be possible.  Note that the only property of $K_m$ used in the
construction above is that $g(X^{(1)}) = g(X) + 1$.
\end{rmk}

As consequences of Theorem~\ref{thm:hopf-annulus}, we obtain the
following.

\begin{cor}
\label{cor:primitive meridian (hopf section)} Let $K_i$, $X_i$,
$K$, $X$ and $X^{(p)}$ be as in the statement of
Theorem~\ref{thm:hopf-annulus}. If $g(X^{(p)}) = g(X)$ for some $p
\ge 1$, and no $X_i$ contains a meridional essential surface $S$
with $\chi(S) \geq 4-2g(X)$, then $K$ admits a primitive meridian.
\end{cor}

\begin{cor}
\label{cor:genus is independant of partition} Let $K_i$, $X_i$,
$X$ and $X^{(p)}$ be as in the statement of
Theorem~\ref{thm:hopf-annulus}. Suppose that no $X_i$ contains a
meridional essential surface $S$ with $\chi(S) \geq
4-2g(X^{(1)})$. Then $g(X^{(1)}) = g(X^{(1)},\del X,\del X^{(1)}
\setminus \del X))$, \ie ,\ there is a minimal genus \hhs\ for
$X^{(1)}$ separating the boundary components.
\end{cor}

\begin{proof}[Proof of Corollary~\ref{cor:primitive meridian (hopf section)}]
By Theorem~\ref{thm:hopf-annulus}, there exists a minimal genus
Heegaard splitting of $X^{(p)}$ (say  $C_1 \cup_{\s^*} C_2$) and a
Hopf-Haken annulus $A^*$ for $\s^*$.  By assumption $g( \s^* ) =
g(X)$. Hence by Proposition~\ref{pro:hopf-haken implies primitive
meridian}, $K$ admits a primitive meridian.
\end{proof}

\begin{proof}[Proof of Corollary~\ref{cor:genus is independant of partition}]
Our assumption implies that conclusion~(2) of
Theorem~\ref{thm:hopf-annulus} cannot hold.  With $p=1$, the
minimal genus Heegaard surface for $X^{(1)}$ described in
conclusion~(1) of Theorem~\ref{thm:hopf-annulus} separates the
boundary components, thus proving the corollary.
\end{proof}

\begin{proof}[Proof of Theorem~\ref{thm:hopf-annulus}]
Recall from the proof of Theorem~\ref{thm:sft}, that we may assume
that each $K_i$ is a prime knot, and both existence and uniqueness
of prime decomposition holds for $K$. We define the complexity of
$X^{(p)}$ to be $(n,p)$ with the lexicographic order (note that
$n$ is the number of prime factors).
Theorem~\ref{thm:hopf-annulus} is proved by induction on this
complexity. As the first step of the induction we prove the
following.

\begin{ass}
\label{ass:1,1}

Suppose that $(n,p) = (1,1)$.  Then Theorem~\ref{thm:hopf-annulus}
holds.
\end{ass}

\begin{proof}[Proof of Assertion~\ref{ass:1,1}]

Let $C_1 \cup_{\s} C_2$ be a minimal genus Heegaard splitting of
$X^{(1)}$ and assume (without loss of generality) that $\del X
\subset \del_- C_1$.  Let $A$ be a Hopf spanning annulus.

\smallskip \noindent
{\bf Case One:  $C_1 \cup_{\s} C_2$ is strongly irreducible.} In
this case we may suppose that each component of $A \cap \s$ is a
\scc\ which is essential in $A$ and $\s$. Minimize $|A \cap \s|$
subject to this constraint (note that $A \cap \s \neq \emptyset$
by (6) of Remarks~\ref{rmk:properties of CBs}).  If $|A \cap \s| =
1$ then we have conclusion (1).  Hence suppose that $|A \cap \s|
\geq 2$.  If some component of $A \cap C_j$ is compressible in
$C_j$ (for $j=1$ or $2$), then by compressing $A$ along that disk
we obtain a disk in $X$ with meridional boundary.  This fact,
together with the irreducibility of $X$ shows that $X$ is a solid
torus, contradicting the assumption of
Theorem~\ref{thm:hopf-annulus} (note that $X_1=X$). Hence every
component of $A \cap C_j$ is incompressible in $C_j$ ($j=1,2$).
This, together with the minimality of $| A \cap \s |$ implies that
each component of $A \cap C_2$ is essential in $C_2$.  Hence there
is a meridian disk of $C_2$ which is disjoint from $A \cap C_2$
((4) of Remarks~\ref{rmk:properties of CBs}).  Then let
$\mathcal{D}$ be a maximal system of mutually disjoint and non
parallel meridian disks of $C_2$ such that $\mathcal{D} \cap (A
\cap C_2) = \emptyset$.  Let $X'$ be the manifold obtained by
cutting $X^{(1)}$ along $A$, and let $\s'$ be the image of $\s$ in
$X'$.  Note that $X'$ is homeomorphic to $X$.  Let $\mathcal{S}'$
be the surface obtained by compressing $\s'$ along (the image of)
$\mathcal{D}$.

\begin{clm}
\label{clm:s' is incompressible} $\mathcal{S}'$ is incompressible
in $X'$.
\end{clm}

\begin{proof}
Suppose that there is a compressing disk $D$ for    $\mathcal{S}'$
in $X'$.  Regard $D$ as a disk in $X^{(1)}$. By the No Nesting
Lemma \cite{no-nestig}, there exists a meridian disk $D'$ in $C_1$
or $C_2$ with $\del D' = \del D$.  By using a standard innermost
disk argument we may suppose that $D' \cap A = \emptyset$. Then,
by isotopy, we may suppose that $D' \cap \mathcal{D} = \emptyset$.
By the maximality of $\mathcal{D}$ we see that $D' \subset C_1$.
Then $D' \cup    \mathcal{D}$ provides a weak reduction for $\s$,
a contradiction.
\end{proof}

If there is a component, say $S'$, of $\mathcal{S}'$ such that
$\del S' \neq \emptyset$ and $S'$ is not boundary parallel in $X'$
then the corresponding surface in $X$ establishes conclusion~(2)
of Theorem~\ref{thm:hopf-annulus} (note that $\chi(S') \geq
4-2g$). Hence we may suppose that each component of $\mathcal{S}'$
is either closed or a boundary parallel annulus. Let $A_1$ and
$A_2$ be the closures of the components of $A \setminus \s$ such
that $A_1 \cap \del X \neq \emptyset$ (hence $A_1 \subset C_1$)
and $A_2$ is adjacent to $A_1$ (\ie ,\ $A_2 \subset C_2$ and $A_2
\cap A_1 \neq \emptyset$).  See Figure~6.2.

\begin{figure}[ht]
\begin{center}
\includegraphics[width=6cm, clip]{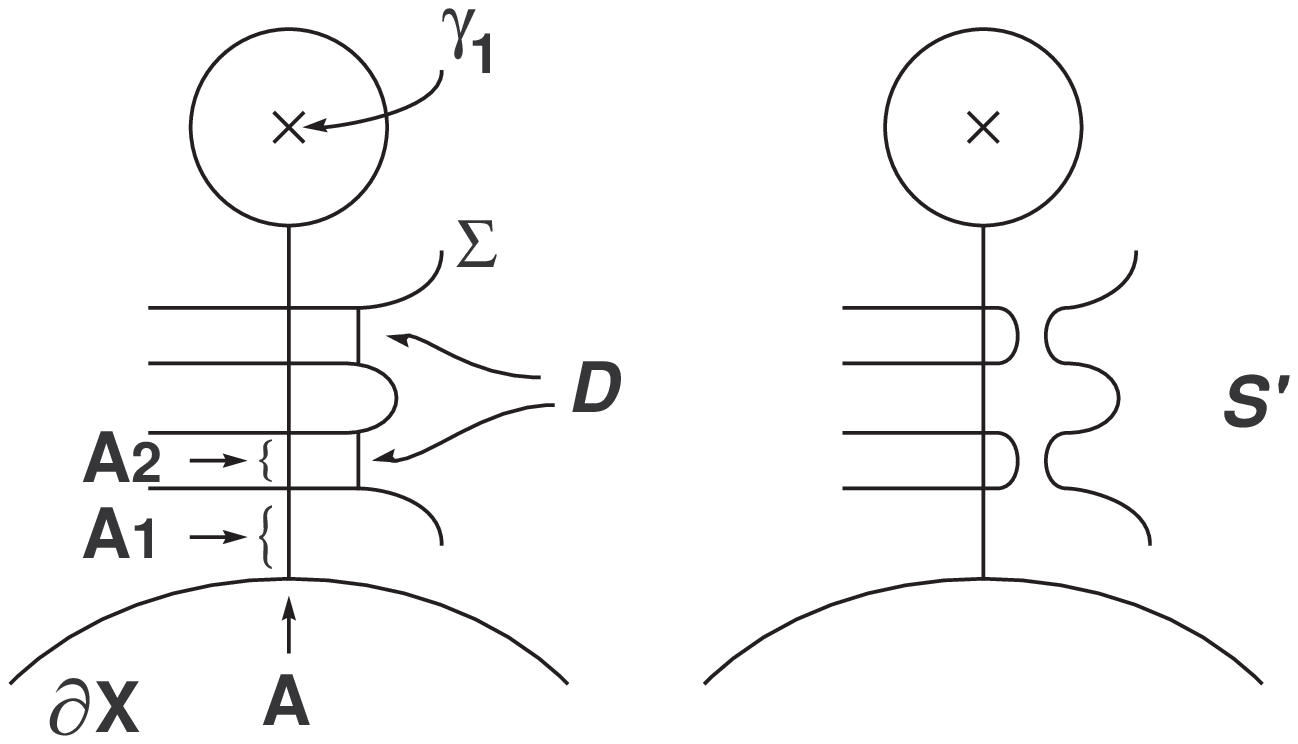}
\end{center}
\begin{center}
Figure 6.2.
\end{center}
\end{figure}

Let $C_2^*$ be the component of $\mbox{cl}(C_2 \setminus
N(\mathcal{D}))$ which contains $A_2$.  Note that $C_2^*$ naturally
inherits a compression body structure from $C_2$
((2) of Remarks~\ref{rmk:properties of CBs}).
Since each component of $\mathcal{S'}$ with non empty
boundary is an annulus, we see that $\del_+ C_2^*$ is a torus;
hence either $C_2^*$ is a solid torus or it is a trivial
compression body homeomorphic to (torus)$\times [0,1]$.

\smallskip \noindent {\bf Subcase 1.a: $C_2^*$ is a solid torus.} See
Figure~6.3.

\begin{figure}[ht]
\begin{center}
\includegraphics[width=8cm, clip]{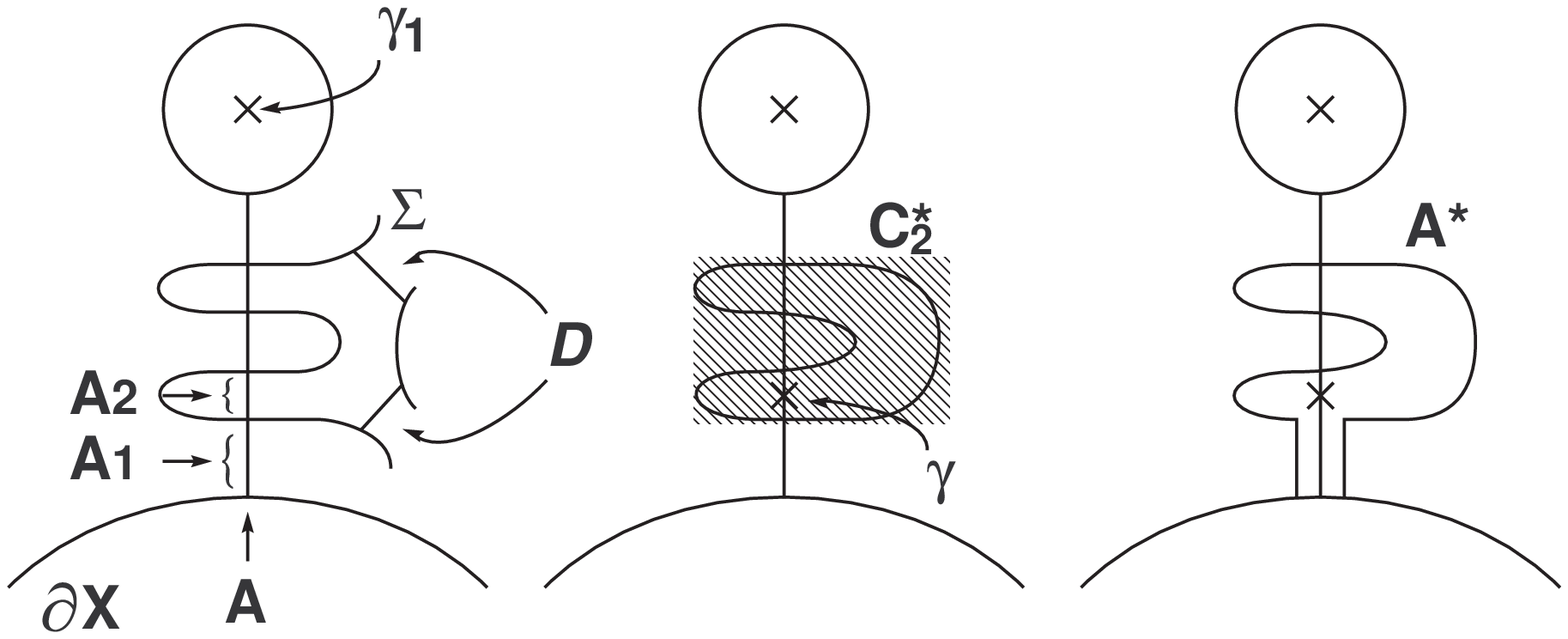}
\end{center}
\begin{center}
Figure 6.3.
\end{center}
\end{figure}

Let $A^*$ be an annulus obtained from $\del C_2^*$ by performing
surgery along $A_1$. We regard $A^*$ as an annulus properly
embedded in $X$ $(=X^{(1)} \cup  N(\gamma_1))$.  If $A^*$ is not
boundary parallel, then $K$ is not prime, contradicting our
assumptions.  Thus $A^*$ is boundary parallel.  If the boundary
parallel region does not contain $C_2^*$ then $X$ is homeomorphic
to the manifold obtained from $C_2^*$ by attaching the parallel
region along an annulus corresponding to $A^*$.  This shows that
$X$ is a solid torus, contradicting the assumptions of
Theorem~\ref{thm:hopf-annulus}.  Hence the parallel region
contains $C_2^*$. This implies that $A_2$ is a longitudinal
annulus in $C_2^*$. Hence the core curve of $A_2$, say $\gamma$,
is a core curve of $C_2^*$. Since $\gamma$ is a core curve of
$C_2^*$ we see that $\s$ is a Heegaard surface of $\mbox{cl}(X
\setminus N(\gamma,C_2^*))$.  Then $A_1$ and a part of $A_2$ form
a Hopf spanning annulus for $\gamma$, and $\s$ intersects this
annulus in a single essential curve. Note that $\gamma_1$ is
isotoped to $\gamma$ along $A$. Hence $\mbox{cl}(X \setminus
N(\gamma,C_2^*))$ is homeomorphic to $X^{(1)}$. These facts
establish conclusion~(1) of Theorem~\ref{thm:hopf-annulus}.

\smallskip \noindent {\bf Subcase 1.b: $C_2^*$ is a trivial
compression body.} Since $\del X \subset C_1$, $\del_- C_2^* =
T_1$ (recall from Section~\ref{sec:sft} that $T_1 = \del
N(\gamma_1))$. Note that there is a vertical annulus $A_2^*$ in
$C_2^*$ such that $A_2^* \cap \del_+ C_2^* = A_1 \cap \del_+
C_2^*$. Then $A_1 \cup A_2^*$ is a Hopf-Haken annulus for $\s$
(Figure~6.4), and we obtain conclusion (1) of
Theorem~\ref{thm:hopf-annulus}.

\begin{figure}[ht]
\begin{center}
\includegraphics[width=6cm, clip]{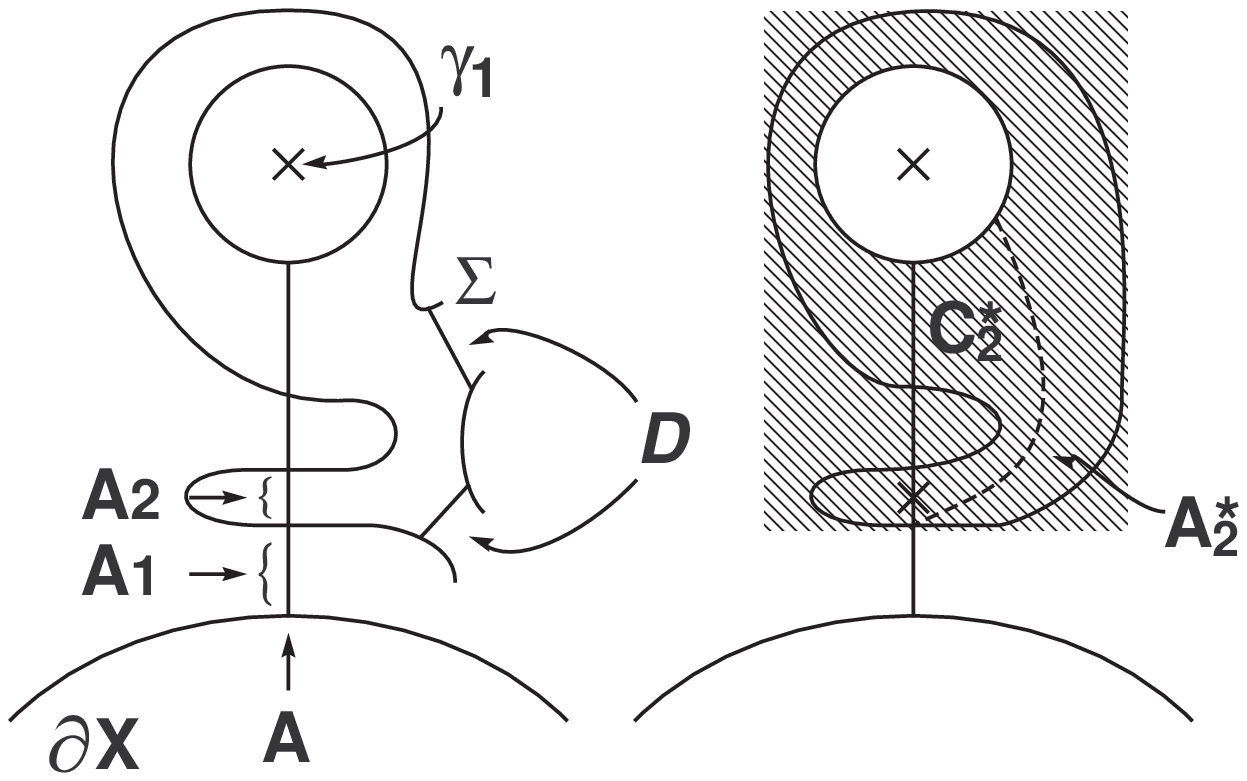}
\end{center}
\begin{center}
Figure 6.4.
\end{center}
\end{figure}

\smallskip \noindent {\bf Case Two: $C_1 \cup_{\s} C_2$ is weakly
reducible.} Let $(C_1^1 \cup_{\s^1} C_2^1) \cup \cdots \cup (C_1^m
\cup_{\s^m} C_2^m)$ be a Scharlemann--Thompson untelescoping of
$C_1 \cup_{\s} C_2$ (Remark~\ref{rmk:scharlemann's ST-untel is
different}). Let $\mathcal{F} = \cup_{i=1}^r F_i$ be the
collection of incompressible surfaces appearing in the
untelescoping, where $F_1,\dots,F_r$ are the connected components.
Note that since $C_1 \cup_{\s} C_2$ is minimal genus, no $F_i$ is
boundary parallel in $X^{(1)}$ (Remark~\ref{rmk:scharlemann's
ST-untel is different}). We may suppose that each component of $A
\cap \mathcal{F}$ is a \scc\ which is essential in $A$ and
$\mathcal{F}$. Minimize $|A \cap \mathcal{F}|$ subject to this
constraint.

\smallskip \noindent {\bf Subcase 2.a: $\mathcal{F} \cap A \neq \emptyset$.}
Without loss of generality we may suppose that $F_1 \cap A \neq
\emptyset$.  Note that $\chi(F_1) \geq 6 -2g$. Let $X'$ be the
manifold obtained by cutting $X^{(1)}$ along $A$ (note that $X'
\cong X$), and $F'$ the image of $F_1$ in $X'$.  Since each
component of $F_1 \cap A$ is essential in $F_1$, $F'$ is
incompressible in $X'$.  Hence if there exists a component of $F'$
which is not boundary parallel in $X'$ then we have conclusion (2)
of Theorem~\ref{thm:hopf-annulus}.

Thus we may suppose that all components of $F'$ are boundary
parallel in $X'$.  In this case, each component of
$F'$ is an annulus.  Let $A'$, $A''$ be the
images of $A$ in $\del X'$, and let $A_X'$, $A_T'$ be the
images of $\del X$, $T_1$ in $\del X'$ respectively.
Note that $A_X'$, $A_T'$ are annuli such that
$\del X' = A_X' \cup A' \cup A_T' \cup A''$.
See Figure~6.5.

\begin{figure}[ht]
\begin{center}
\includegraphics[width=2.5cm, clip]{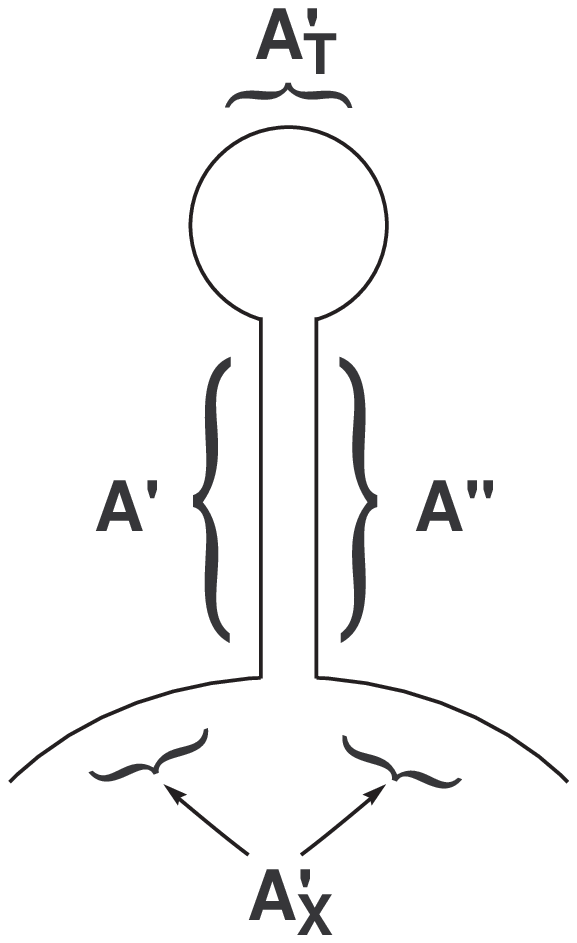}
\end{center}
\begin{center}
Figure 6.5.
\end{center}
\end{figure}

Suppose that there exists a component of $F'$, say $A^*$, such
that one component of $\del A^*$ is contained in $A'$, and the
other in $A''$.   Then $A^*$ is parallel to $A_X'$ or $A_T'$.  We
may suppose that $A^*$ is the outermost component with that
property.

Suppose that the components of $\del A^*$ are not identified in
$F_1$ (see Figure~6.6 (i)). Let $A^{**}$ be the component of $F'$
which is adjacent to $A^*$ and is contained in the parallel region
bounded by $A^*$. Then it is easy to see that $A^{**}$ is parallel
to an annulus in $A'$ or $A''$. However this contradicts the
minimality of $|A \cap \mathcal{F}|$. Hence the components of
$\del A^*$ are identified in $F_1$. This shows that $F_1$ is a
boundary parallel torus in $X^{(1)}$, contradicting
Remark~\ref{rmk:scharlemann's ST-untel is different}.

\begin{figure}[ht]
\begin{center}
\includegraphics[width=6cm, clip]{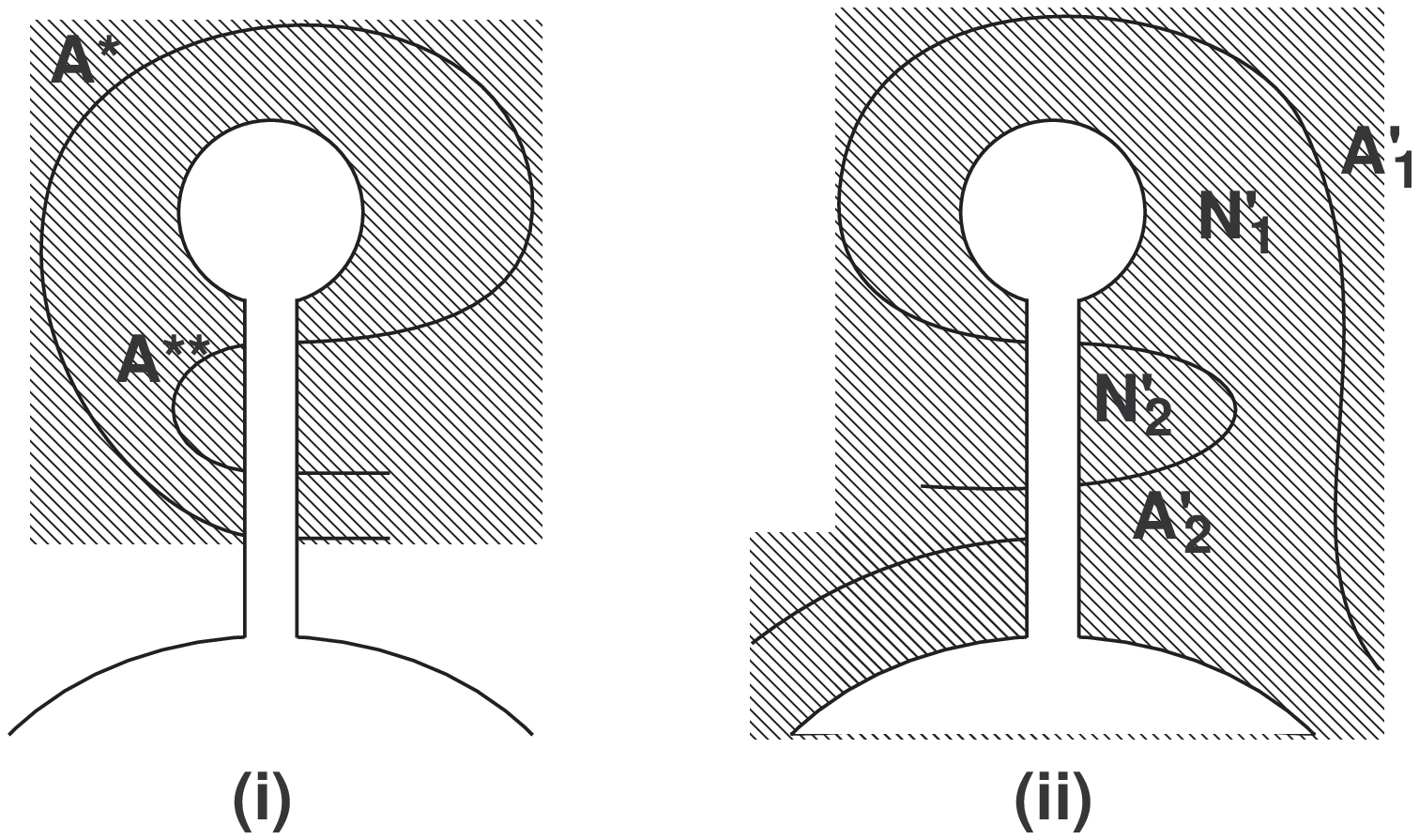}
\end{center}
\begin{center}
Figure 6.6.
\end{center}
\end{figure}

Now suppose that for each component $A^*$ of $F'$ we have that
either $\del A^* \subset A'$ or $\del A^* \subset A''$.  Let us
fix $A_1'$, a component of $F'$.  We may suppose (without loss of
generality) that $\del A_1' \subset A'$.  Let $A_2'$ be another
component of $F'$ such that $A_2'$ is adjacent to  $A_1'$ in
$F_1$.  Then $\del A_2' \subset A''$.  Let $N_1'$ be the parallel
region between $A_1'$ and a subsurface of $\del X'$.  By the
minimality of $|\mathcal{F} \cap A|$ we see that $N_1' \cap \del
X'$ is not contained in $A'$. Hence $N_1' \cap \del X'$ contains
$A''$, and $A_2' \subset N_1'$.  See Figure~6.6 (ii). However,
this implies that $A_2'$ admits a parallel region $N_2'$ such that
$N_2' \cap \del X'$ is contained in $A''$, contradicting the
minimality of $|\mathcal{F} \cap A|$. Hence there exists a
component of $F'$ which is not boundary parallel and this
establishes conclusion~(2) of Theorem~\ref{thm:hopf-annulus}.

\smallskip \noindent {\bf Subcase 2.b: $\mathcal{F} \cap A = \emptyset$.} We
may assume that $\del X \subset C_1^1$.  In this case, we
basically apply the arguments of Case One to the strongly
irreducible Heegaard splitting $C_1^1 \cup_{\s^1} C_2^1$.  Let $M
= C_1^1 \cup C_2^1$ (\ie ,\ the ambient manifold of the Heegaard
splitting $C_1^1 \cup_{\s^1} C_2^1$).  Then $A \subset M$.  Since
$C_1^1 \cup_{\s^1} C_2^1$ is strongly irreducible, we may suppose
that each component of $A \cap \s^{1}$ is a \scc\ which is
essential in both $A$ and $\s^{1}$. Then by Remark
\ref{rmk:properties of CBs}(6), $A \cap \s^{1} \neq \emptyset$ and
each component of $A \cap C_{i}^1$ is incompressible in $C_i^1$
($i=1,2$).  If $|A \cap \s^{1}|= 1$, then by
Proposition~\ref{pro:haken annulus after amalgamation} we have
conclusion (1) of Theorem ~\ref{thm:hopf-annulus}.  Hence we may
suppose that $|A \cap\s^{1}| \geq 2$.

Let $A_1$, $A_2$ be the closures of the components of $A \setminus
\s^1$ such that $A_1 \cap \del X \neq \emptyset$ (hence $A_1
\subset C_1^1$) and $A_2 \subset C_2^1$ with $A_2 \cap A_1 \neq
\emptyset$.  Then by arguments of Case One we can show that
there exists a union of mutually disjoint meridian disks of
$C_2^1$ (say $\mathcal{D}$) which is non empty and $\mathcal{D}
\cap (A \cap C_2^1) = \emptyset$.  Let $M'$ be the manifold
obtained from $M$ by cutting along $A$, and let $\s'$ be the image
of $\s^1$ in $M'$. Let $\mathcal{S}'$ be the surface obtained by
compressing $\s'$ along $\mathcal{D}$. Then by using the argument
of the proof of Claim ~\ref{clm:s' is incompressible} in Case One
we can show that $\mathcal{S}'$ is incompressible in $M'$.

If there is a component of $\mathcal{S'}$ which is neither closed
nor boundary parallel in $M'$ then it is easy to see that the
surface establishes conclusion~(2) of
Theorem~\ref{thm:hopf-annulus}. Therefore we may suppose that each
component of $\mathcal{S'}$ is either closed or boundary parallel
in $M'$.

Let $C_2^*$ be the component of $\mbox{cl}(C^1_2 \setminus
N(\mathcal{D}))$ which contains $A_2$. Then $C_2^*$ is a
compression body and $\del_+ C_2^*$ is a torus (for a detailed
argument, see Case One.)  Hence, we have the following three
subcases:

\begin{description}
\item[Subcase A] $C_2^*$ is a solid torus,
\item[Subcase B] $C_2^*$ is a trivial compression body
homeomorphic to
(torus)$\times [0,1]$ with $\del_- C_2^* = T_1$, or---
\item[Subcase C] $C_2^*$ is a trivial compression body
homeomorphic to
(torus)$\times [0,1]$ with $\del_- C_2^* = F_k \subset
\mathcal{F}$.
\end{description}

For Subcases A and B, we apply the argument of Subcases 1.a and
1.b in Case One to show that there is a Hopf-Haken annulus for
$\s^1$ in $M$ with, possibly, the position of $\gamma_1$ changed.
Then by Proposition~\ref{pro:haken annulus after amalgamation}, we
have conclusion (1) of Theorem~\ref{thm:hopf-annulus}. We may
therefore suppose that Subcase C holds. See Figure~6.7.

\begin{figure}[ht]
\begin{center}
\includegraphics[width=7cm, clip]{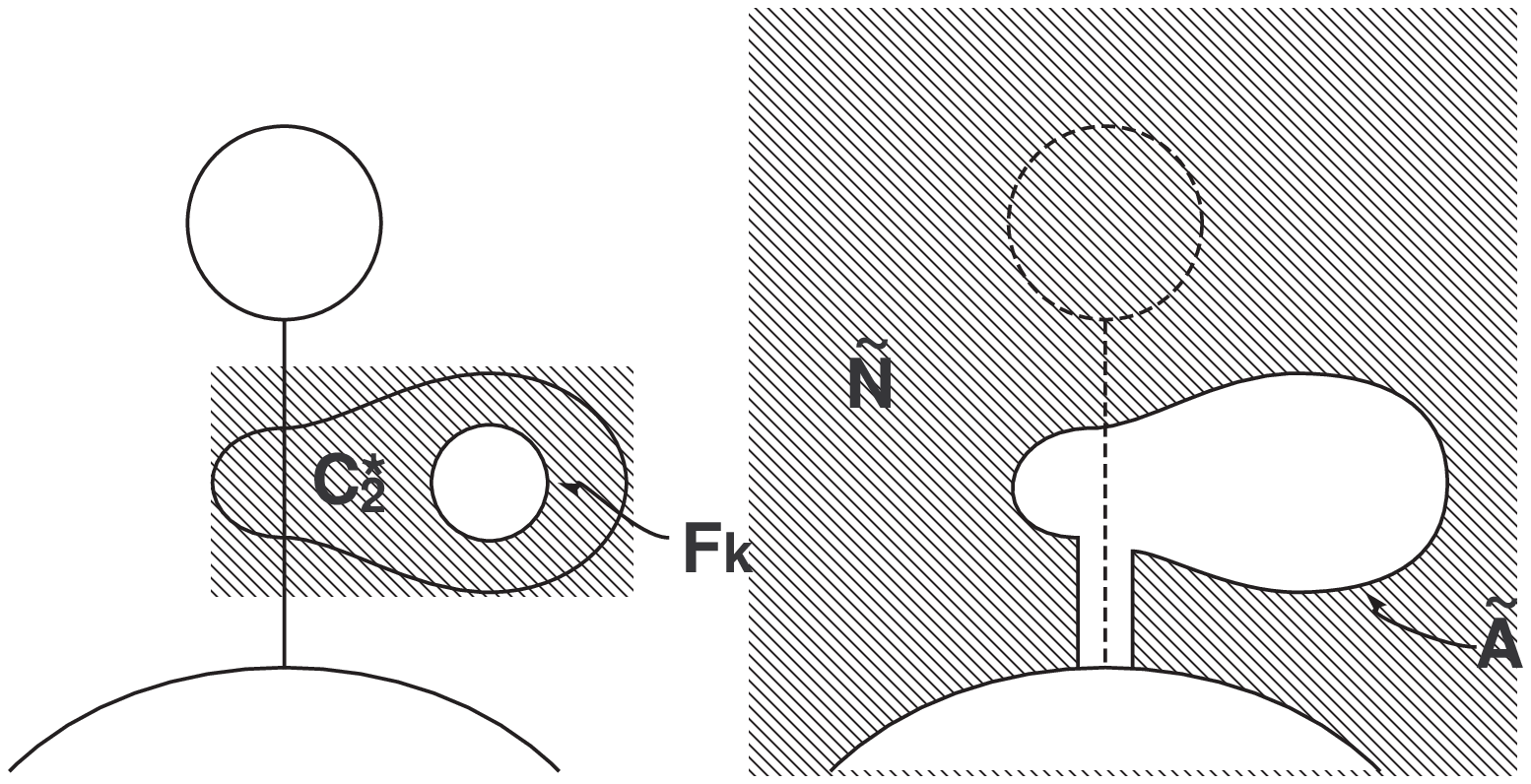}
\end{center}
\begin{center}
Figure 6.7.
\end{center}
\end{figure}

Let $\widetilde A$ be the annulus obtained from $\del_+ C_2^*$ by
performing surgery along $A_1$.  Regard $\widetilde A$ as an
annulus in $X (=X^{(1)} \cup N(\gamma_1))$.  If $\widetilde A$ is
not boundary parallel then we have a contradiction to $K$ being
prime. Thus $\widetilde A$ is boundary parallel. Let $\widetilde
N$ be the parallel region. Then $\widetilde N$ is a solid torus
and $\widetilde N \cap \del X = \del \widetilde N \cap \del X$ is
a longitudinal annulus, that is, the annulus wraps around $\widetilde{N}$
once.

Suppose that $\widetilde N \supset F_k$. Since $\del_+ C_2^*$ and
$F_k$ are parallel, this implies that $F_k$ bounds a solid torus
in $X$ which is isotopic to $\widetilde N$. Since $F_k$ is
incompressible in $X^{(1)}$, $\gamma_1$ must be contained in this
solid torus.  Note that $A$ connects $\gamma_1$ to $\del X$.
However this implies that $A \cap F_k \neq \emptyset$,
contradicting the condition of Subcase 2.b. Hence we have
$\widetilde N \cap F_k = \emptyset$ (see the right side of
Figure~6.7). This shows that $\del_+ C_2^*$ (hence also $F_k$
which is parallel to it) is boundary parallel in $X$. 
By Remark~\ref{rmk:scharlemann's ST-untel is different}, $F_k$ is
not boundary parallel in $X^{(1)}$, therefore $\gamma_1$ must be
contained in the parallel region between $F_k$ and $\del X$. Then
recall  that $A \cap F_k = \emptyset$ (the condition of Subcase
2.b), and $\partial A$ consists of $\gamma_1$ and a meridian curve
of $\partial X$. This shows that $M \cong D(2)$ (with $D(2) \cong$
(disk with two holes) $\times S^1$, as defined in the end of
Section~\ref{sec:background}).

As conclusions of the above we have the following:

\begin{enumerate}
\item $F_k$ is a boundary parallel torus in $X$, and---
\item $M$ is homeomorphic to a (disk with two holes)
$\times S^1$, where $\del M = \del X \cup T_1 \cup F_k$.
\end{enumerate}

Let $X^* = \mbox{cl}(X^{(1)} \setminus M)$ $(\cong X)$. Note that
$M \cap X^* = F_k$. By Proposition~\ref{pro:untel to a conn sep
surfce implies weak reduction}, there is a weakly reducing
collection of disks $\Delta$ for $\s$ such that $\widehat \s (
\Delta ) = F_k$. Then by Proposition~\ref{prop:amalgamation of
minimal genus: connected case}, we see that $g(X^{(1)}) =
g(M)+g(X^*)-1$. By (2) of 
Remarks~\ref{rmk:D(n) admits a haken annulus},
there exists a minimal genus Heegaard splitting $\bar C'_1
\cup_{\bar \s'} \bar C'_2$ of $M$ such that $\bar \s' \cap A$
consists of a single \scc\ that is essential in $A$. Take the
amalgamation of $\bar C'_1 \cup_{\bar \s'} \bar C'_2$ and a
minimal genus Heegaard splitting of $X^*$. Then by 
Lemma~\ref{lem:genus after amalgamation}, and 
Proposition~\ref{pro:haken annulus after amalgamation}, we see that
conclusion (1) of Theorem~\ref{thm:hopf-annulus} holds.

    This completes the proof of Assertion~\ref{ass:1,1}.
\end{proof}

\begin{ass}
\label{ass:n=1}

Suppose that $n=1$.  Then Theorem~\ref{thm:hopf-annulus} holds.
\end{ass}

\begin{proof}[Proof of Assertion~\ref{ass:n=1}]

Assertion~\ref{ass:n=1} is proved by inducting on $p$.  Assertion
\ref{ass:1,1} gives the first step of the induction.  Hence we
suppose that $p \geq 2$.  Let $C_1 \cup_{\s^{(p)}} C_2$ be a
minimal genus Heegaard splitting for $X^{(p)}$ and assume that
$\del X \subset \del_- C_1$.  Recall (from Section~\ref{sec:sft})
that $\del X^{(p)} = \del X \cup T_1 \cup \cdots \cup T_p$.  Let
$A^1,\dots,A^p$ be mutually disjoint Hopf spanning annuli between
$\del X$ and $T_1,\dots,T_p$, respectively, and $A = \cup_{j=1}^p
A^j$. See Figure~6.8.
By (6) of Remarks~\ref{rmk:properties of CBs}, we see that for each
$j,$ $A^j \cap \s^{(p)} \neq \emptyset$.

\begin{figure}[ht]
\begin{center}
\includegraphics[width=6cm, clip]{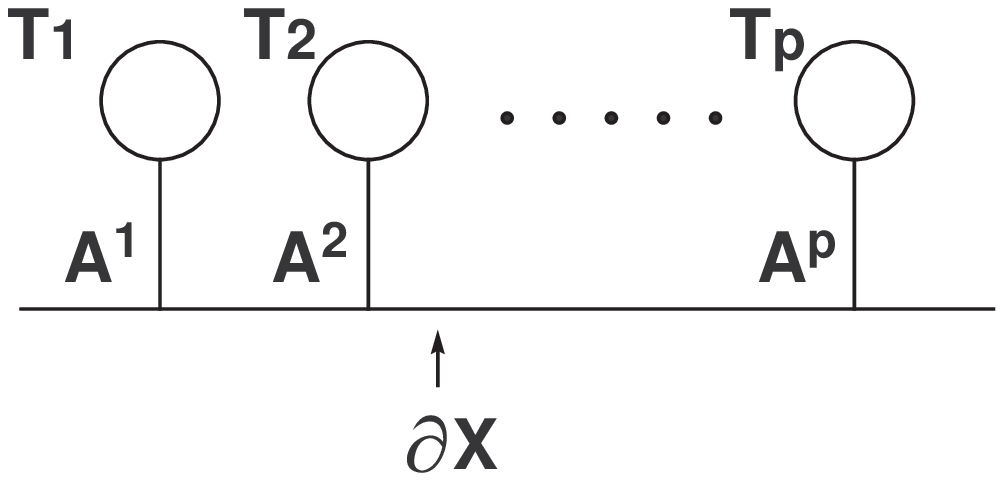}
\end{center}
\begin{center}
Figure 6.8.
\end{center}
\end{figure}

\smallskip \noindent
{\bf Case One: $C_1 \cup_{\s^{(p)}} C_2$ is strongly irreducible.} 
In this case we may suppose that $A \cap \s^{(p)}$ consists of
\scc s that are essential in $A$ and $\s^{(p)}$.  Subject to that
constraint we minimize $|A \cap \s^{(p)}|$. If $| A^j \cap
\s^{(p)} |=1$ for some $j$, then we have conclusion (1) of
Theorem~\ref{thm:hopf-annulus}. Hence suppose that for each $j$,
$|A^j \cap \s^{(p)}| \geq 2.$  The following argument is quite
similar to that of Step One in the proof of
Assertion~\ref{ass:1,1}; hence we will only sketch it here.

Let $\mathcal{D}$ be a maximal collection of meridian disks for
$C_2$ so that $\mathcal{D} \cap (A \cap C_2) = \emptyset$; by (4)
of Remarks~\ref{rmk:properties of CBs}, $\mathcal{D} \neq
\emptyset$. Let $X'$ be the manifold obtained by cutting $X^{(p)}$
along $A$, and let $\s'$ be the image of $\s^{(p)}$ in $X'$.  Let
$\mathcal{S}'$ be the surface obtained from $\s'$ by compressing
along $\mathcal{D}$.  Then by the same argument as in the proof of
Claim~\ref{clm:s' is incompressible} of the proof of Assertion
\ref{ass:1,1}, we see that $\mathcal{S}'$ is incompressible in
$X'$.

If there exists a component of $\mathcal{S}'$ with non-empty
boundary that is not boundary parallel in $X'$ then we have
conclusion (2) of Theorem~\ref{thm:hopf-annulus}.  Thus we may
assume that each component of $\mathcal{S}'$ is either closed or
boundary parallel.  (Hence each component of $\mathcal{S}'$ with
non-empty boundary is a boundary parallel annulus.)  Let $A_1$,
$A_2$ be the closures of the components of $A^1 \setminus
\s^{(p)}$ such that $A_1 \cap \del X \neq \emptyset$ (hence $A_1
\subset C_1$) and that $A_2 \subset C_2$ with $A_2 \cap A_1 \neq
\emptyset$.  Let $C_2^*$ be the component of $\mbox{cl}(C_2
\setminus N(\mathcal{D}))$ which contains $A_2.$  Then $C_2^*$ is
either a solid torus or homeomorphic to (torus)$\times [0,1]$.

\smallskip \noindent
{\bf Subcase 1.a: $C_2^*$ is a solid torus.}   In this subcase the
argument of Subcase 1.a of the proof of Assertion~\ref{ass:1,1}
can be applied to show that we can isotope $\gamma_1$ along $A^1$
to a core curve of $A_2$ to show that there exists a minimal genus
Heegaard splitting of $X^{(p)}$ that intersects a Hopf spanning
annulus in a single essential \scc\ (conclusion (1) of
Theorem~\ref{thm:hopf-annulus}).

\smallskip \noindent
{\bf Subcase 1.b: $C_2^*$ is homeomorphic to (torus)$\times 
[0,1].$}   In this case the argument of Subcase 1.b of the proof of
Assertion~\ref{ass:1,1} applies to show that conclusion~(1) of
Theorem~\ref{thm:hopf-annulus} holds.

\smallskip \noindent
{\bf Case Two: $C_1 \cup_{\s^{(p)}} C_2$ is weakly reducible.} 
Let

$$(C_1^1 \cup_{\s^1} C_2^1)\cup\cdots\cup (C_1^m \cup_{\s^m} C_2^m)$$

\noindent be a Scharlemann-Thompson untelescoping of $C_1
\cup_{\s^{(p)}} C_2$, where $\del X \subset \del_- C_1^1$.  Let
$\mathcal{F}$ be the union of the essential surface appearing in
the above untelescoping (Remark~\ref{rmk:scharlemann's ST-untel is
different}) and $\{F_i\}_{i=1}^r$ its connected components. We may
suppose that each component of $\mathcal{F} \cap A$ is a \scc\
that is essential in both $\mathcal{F}$ and $A$, and $|\mathcal{F}
\cap A|$ is minimal subject to this constraint.  Then we have the
following subcases.

\smallskip \noindent
{\bf Subcase 2.a: $\mathcal{F} \cap A \neq \emptyset.$} Without
loss of generality we may suppose that $F_1 \cap A^1 \neq
\emptyset$.  Let $F'$ be the image of $F_1$ in $X'$ ($X'$ as in
Case One). Since each component of $F_1 \cap A$ is essential in
$F_1$, $F'$ is incompressible in $X'$.  Hence if there exists a
component of $F'$ which is not boundary parallel, then we have
conclusion (2) of Theorem~\ref{thm:hopf-annulus}.  We suppose that
each component of $F'$ is boundary parallel in $X'$ (see
Figure~6.9).

\begin{figure}[ht]
\begin{center}
\includegraphics[width=8cm, clip]{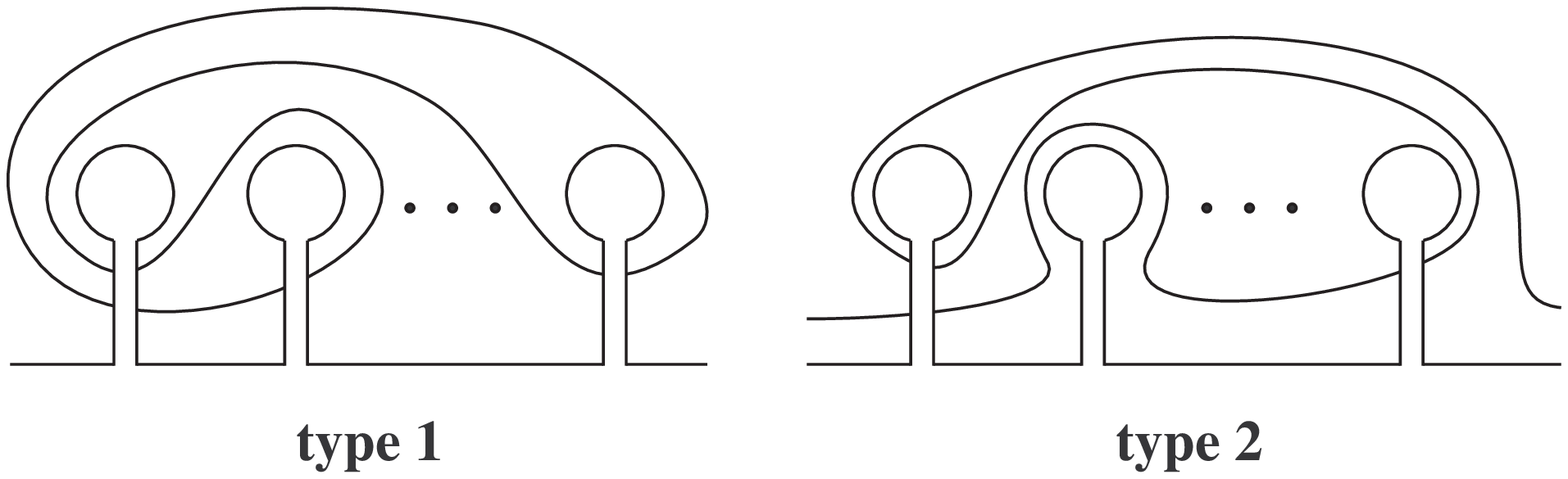}
\end{center}
\begin{center}
Figure 6.9.
\end{center}
\end{figure}

This implies that $F'$ can be isotoped to be contained in $N(\del
X',X')$.  Here we note that $X^{(p)}$ admits the decomposition

$$X^{(p)} = X^* \cup_T D(p+1),$$

\noindent where $X^* \cong X$, $D(p+1)$ is a (disk with $p+1$
holes)$\times S^1$ (see Section~\ref{sec:background}), and a
meridian of $X^*$ is identified with $\{*\} \times S^1$, where $*$
is a point on the boundary of the disk.

Since $F'$ is contained in $N(\del X',X')$, we may suppose that
$F_1$ is contained in $\mbox{int}(D(p+1))$.  Then by
\cite[VI.34,VI.18]{jaco}, we see that 
$F_1$ is vertical, \ie ,\ 
$F_1$ is isotopic to a torus
of the form (\scc) $\times S^1$ in $D(p+1)$.  Hence $F_1$ gives
the following decomposition:

$$X^{(p)} = X^{(p_1)} \cup_{F_1} D(p_2),$$

\noindent where $p_1+p_2=p+1$. We say that $F_1$ is of type 1 if
$\del X \subset \del X^{(p_1)}$ and of type 2 if $\del X \subset
\del D(p_2)$.  See Figure~6.9.

Suppose that $F_1$ is of type 1. Here we note that $p_2 \ge 2$
(Remark~\ref{rmk:scharlemann's ST-untel is different}); hence $p_1
< p$ and we may apply the inductive hypothesis on $X^{(p_1)}$,
that is, there is a minimal genus Heegaard splitting $C_1''
\cup_{\s''} C_2''$ of $X^{(p_1)}$ and a Hopf spanning annulus
$A_1'' \subset X^{(p_1)}$ between $\del X$ and $F_1$ so that $\s''
\cap A_1''$ is a single \scc\ which is essential in $A_1''$. On
the other hand, 
by (2) of Remarks~\ref{rmk:D(n) admits a haken annulus}
there is a minimal genus Heegaard splitting of $D(p_2)$ and an
essential annulus $A_2'' \subset D(p_2)$ such that the
intersection of the Heegaard surface and $A_2''$ is a single \scc\
which is essential in $A_2''$, and $A_2'' \cap F_1 = A_1'' \cap
F_1$. Then by applying Lemma~\ref{lem:genus after amalgamation},
Propositions~\ref {prop:amalgamation of minimal genus: connected
case} and~\ref{pro:combined haken annuli} to these Heegaard
splittings we see that conclusion (1) of
Theorem~\ref{thm:hopf-annulus} holds.

Suppose that $F_1$ is of type~2. By Remark~\ref{rmk:scharlemann's
ST-untel is different}, we see that $\del D(p_2)$ contains at
least one $T_i$. By changing subscripts if necessary, we may
suppose that $T_1 \subset \del D(p_2)$. We take a minimal genus
Heegaard splitting of $D(p_2)$ and an essential annulus $A_2''
\subset D(p_2)$ such that $A_2''$ intersects the Heegaard surface
in a single simple closed curve which is essential in $A_2''$ and
that $A_2''$ connects $\del X$ to $T_1$ (Remarks~\ref{rmk:D(n)
admits a haken annulus}(2)). Then take an amalgamation of this
Heegaard splitting and a minimal genus Heegaard splitting of
$X^{(p_1)}$. By Lemma~\ref{lem:genus after amalgamation},
Propositions~\ref {prop:amalgamation of minimal genus: connected
case} and~\ref{pro:haken annulus after amalgamation}, we see that
we have conclusion (1) of Theorem~\ref{thm:hopf-annulus}.

\smallskip \noindent
{\bf Subcase 2.b: $\mathcal{F} \cap A = \emptyset$.} We apply 
arguments similar to that of Case One to the strongly irreducible
Heegaard splitting $C_1^1 \cup_{\s^1} C_2^1$.  Let $M = C_1^1 \cup
C_2^1$. By the argument of Case One of the proof of
Assertion~\ref{ass:1,1}, we see that after isotopy each component
$A \cap C_2^1$ is essential in $C_2^1$.   Let $\mathcal{D}$ be a
maximal collection of meridian disks for $C_2^1$ such that
$\mathcal{D} \cap (A \cap C_2^1) = \emptyset$ (as above,
$\mathcal{D} \neq \emptyset$). Let $M'$ be the manifold obtained
by cutting $M$ along $A$, and let $\s'$ be the image of $\s^1$ in
$M'$. Let $\mathcal{S}'$ be the surface obtained by compressing
$\s'$ along $\mathcal{D}$. Then by using the argument of
Claim~\ref{clm:s' is incompressible} in the proof of
Assertion~\ref{ass:1,1}, we can show that $\mathcal{S}'$ is
incompressible in $M'$.

If there is a component of $\mathcal{S}'$ with non-empty boundary
that is not boundary parallel, then we have conclusion (2) of
Theorem~\ref{thm:hopf-annulus}.   Suppose then that each component
of $\mathcal{S}'$ is either closed or boundary parallel in $X'$.
As before, let $A_1$, $A_2$ be the closures of components of $A^1
\setminus \s^1$ such that $A_1 \cap \del X \neq \emptyset$, and
that $A_2 \subset C_2^1$ with $A_2 \cap A_1 \neq \emptyset$.  Let
$C_2^*$ be the component of $\mbox{cl}(C_2^1 \setminus
N(\mathcal{D}))$ containing $A_2$. As in the proof of
Assertion~\ref{ass:1,1}, we have the following subcases:

\smallskip \noindent
{\bf  Subcase A:}  $C_2^*$ is a solid torus.

\smallskip \noindent
{\bf  Subcase B:}  $C_2^*$ is a trivial compression body
homeomorphic to (torus)$\times [0,1]$ with $\del_- C_2^* = T_i$.

\smallskip \noindent
{\bf  Subcase C:} $C_2^*$ is a trivial compression body
homeomorphic to (torus)$\times [0,1]$ with $\del_- C_2^* = F_k
\subset \mathcal{F}$.

For Subcases A, B we apply the arguments of Subcases~1.a and 1.b
in Case One, using Proposition~\ref{pro:haken annulus
after amalgamation} conclude that we have conclusion (1) of
Theorem~\ref{thm:hopf-annulus}.

Suppose Subcase C holds (see Figure~6.7). Let $\widetilde A$ be an
annulus obtained from $\del_+ C_2^*$ by performing surgery along
$A_1$.  Regard $\widetilde A$ as an annulus in $X$ ($= X^{(p)}
\cup (\cup_{i=1}^p N(\gamma_i))$).  If $\widetilde A$ is not
boundary parallel, then we have a contradiction to our assumption
that $n=1$. Therefore we may suppose that $\widetilde A$ is
boundary parallel in $X$. Then by using the argument of Subcase~C of
Subcase~2.b of the proof of Assertion~\ref{ass:1,1}, we see that
the following hold.  (Note that $A \subset M \cup (\cup_{i=1}^p
N(\gamma_i))$.)

\begin{enumerate}
    \item $F_k$ is boundary parallel in $X$, and---
    \item $M$ is homeomorphic to
    $D(p+1)$ where $\del M = \del X \cup
    F_k \cup (T_1 \cup \cdots \cup T_p)$.
\end{enumerate}

\noindent This gives the following decomposition:

$$X^{(p)} = X \cup_{F_k} D(p+1).$$

\noindent Then by the inductive hypothesis, Lemma~\ref{lem:genus
after amalgamation}, Remark~\ref{rmk:D(n) admits a haken annulus},
Propositions~\ref{prop:amalgamation of minimal genus: connected
case}, and~\ref{pro:haken annulus after amalgamation}, we see that
conclusion (1) of Theorem~\ref{thm:hopf-annulus} holds (for
detailed arguments, see the last paragraph of Subcase~2.a of this
proof).

This completes the proof of Assertion~\ref{ass:n=1}.
\end{proof}

We now complete the proof of Theorem~\ref{thm:hopf-annulus}.
Recall that the proof is carried out by induction on $(n,p)$ with
the lexicographic ordering.  By Assertion~\ref{ass:n=1} we may
assume that $n \geq 2$.  Let $C_1 \cup_{\s} C_2$ be a minimal
genus Heegaard splitting of $X^{(p)}$.  By Theorem~\ref{thm:sft}
we either have conclusion (2) of Theorem~\ref{thm:hopf-annulus},
or $C_1 \cup_{\s} C_2$ weakly reduces to a \sft\ $T$, which gives
the following decomposition of $X^{(p)}$:

$$X^{(p)} = \mathcal{X}_1^{(p_1+1)} \cup_T \mathcal{X}_2^{(p_2)},$$

\noindent where $\del X \subset \del \mathcal{X}_1^{(p_1+1)}$,
$p_1+p_2 = p$ and $\mathcal{X}_1 = E(\mathcal{K}_1)$,
$\mathcal{X}_2 = E(\mathcal{K}_2)$ (possibly $\mathcal{K}_1$ or
$\mathcal{K}_2$ being the trivial knot in $S^3$) and $K =
\mathcal{K}_1 \# \mathcal{K}_2$.

\smallskip \noindent
{\bf Case 1:} $\mathcal{K}_1$ is the trivial knot in $S^3$. In
this case $\mathcal{X}_1^{(p_1+1)}$ is (disk with $p_1+1$ holes)
$\times S^1$, $p_1 \geq 1$ (Remark \ref{rmk:scharlemann's ST-untel
is different}).  By changing the subscripts if necessary, we may
suppose that $T_1 \subset \del \mathcal{X}_1^{(p_1+1)}$. 
By (2) of 
Remarks~\ref{rmk:D(n) admits a haken annulus} there exists a
minimal genus Heegaard splitting $C_1^1 \cup_{\s^1} C_2^1$ of
$\mathcal{X}_1^{(p_1+1)}$ and an essential annulus $A$ in
$\mathcal{X}_1^{(p_1+1)}$ such that $A$ connects $\del X$ to
$T_1$, $A \cap \del X$ is a meridian, and $\s^1 \cap A$ is a
single \scc\ that is essential in $A$.  Then by
Lemma~\ref{lem:genus after amalgamation},
Propositions~\ref{prop:amalgamation of minimal genus: connected
case}, and~\ref{pro:haken annulus after amalgamation}, we see that
conclusion (1) of Theorem~\ref{thm:hopf-annulus} holds.

\smallskip \noindent
{\bf Case 2:} $\mathcal{K}_2$ is the trivial knot in $S^3$. In
this case $\mathcal{X}_2^{(p_2)}$ is (disk with $p_2$
holes)$\times S^1$, $p_2 \geq 2$
(Remark~\ref{rmk:scharlemann's ST-untel is different}). 
By changing the subscripts if
necessary, we may suppose that $T_1 \subset \del
\mathcal{X}_2^{(p_2)}$. Since $p_1 \geq 2$, we have that $p_1 \leq p-1$, and
this allows us to apply the inductive hypothesis to
$\mathcal{X}_1^{(p_1+1)}$.  
If conclusion~(2) holds for $\mathcal{X}_1^{(p_1+1)}$,
then conclusion (2) holds for $X^{(p)}$ 
(note that $K=\mathcal{K}_1$, and
$g(\mathcal{X}_1^{(p_1+1)}) \le g(X^{(p)})$). Hence we may suppose that
conclusion (1) holds for $\mathcal{X}_1^{(p_1+1)}$, 
{\it i.e.}, there exists
a minimal genus Heegaard splitting (say $C_1^1 \cup_{\s^1} C_2^1$)
of $\mathcal{X}_1^{(p_1+1)}$ admitting 
a Hopf--Haken annulus $A^*$, so that $A^*
\cap T \neq \emptyset$.  By (2) of
Remarks~\ref{rmk:D(n) admits a haken annulus} there exists a
minimal genus Heegaard splitting $C_1^2 \cup_{\s^2} C_2^2$ of
$\mathcal{X}_2^{(p_2)}$, and an essential annulus $A$ in
$\mathcal{X}_2^{(p_2)}$ such that $A$ connects $T_1$ to $T$, 
$A \cap T = A^* \cap T$, 
and $A\cap \s^2$ is
a \scc\ which is essential in $A$. Then by  Lemma~\ref{lem:genus
after amalgamation} and Propositions~\ref{prop:amalgamation of
minimal genus: connected case} and~\ref{pro:combined haken
annuli}, we see that $C_1^1 \cup_{\s^1} C_2^1$ and $C_1^2
\cup_{\s^2} C_2^2$ are amalgamated
to give a 
minimal genus 
Heegaard splitting of $X^{(p)}$ admitting a Hopf--Haken
annulus; hence conclusion (1) of Theorem~\ref{thm:hopf-annulus}
holds.

\smallskip \noindent
{\bf Case 3:} Both $\mathcal{K}_1$ and $\mathcal{K}_2$ are
non-trivial knots.  Then we have the following subcases:

\smallskip \noindent
{\bf Subcase 3.a: $p_2 = 0$.}  Since $p_1 + p_2 = p > 0$, this
implies that $p_1 > 0$. Note that $T_1$ is contained in
$\mathcal{X}_1^{(p_1+1)}$. Recall that we used lexicographic order
on $(n,p)$ and since $\mathcal{K}_2$ is a non-trivial knot, the
number of prime factors of $\mathcal{K}_1$ is strictly less than
that of $K$. We may therefore apply the inductive hypothesis to
$\mathcal{X}_1^{(p_1 + 1)}$.
If conclusion (2) of Theorem~\ref{thm:hopf-annulus} 
holds for $\mathcal{X}_1^{(p_1 + 1)}$, then 
conclusion (2) holds for $X^{(p)}$. 
Suppose that conclusion (1) holds for 
$\mathcal{X}_1^{(p_1 + 1)}$.
Then we obtain a minimal genus Heegaard
splitting $C_1^1 \cup_{\s^1} C_2^1$ of $\mathcal{X}_1^{(p_1 + 1)}$
and an essential annulus $A$ connecting $\del X$ to $T_1$ such
that $A \cap \s^1$ is a single simple closed curve which is
essential in $A$. Then by Lemma~\ref{lem:genus after
amalgamation} and Propositions~\ref{prop:amalgamation of minimal
genus: connected case} and~\ref{pro:haken annulus after
amalgamation}, we see that conclusion (1) of
Theorem~\ref{thm:hopf-annulus} holds.

\smallskip \noindent
{\bf Subcase 3.b: $p_2 > 0$.}  We may assume that $T_1$ is
contained in $\mathcal{X}_2^{(p_2)}$.  Apply the inductive
hypothesis to both $\mathcal{X}_1^{(p_1+1)}$ and
$\mathcal{X}_2^{(p_2)}$.  If conclusion~(2) of
Theorem~\ref{thm:hopf-annulus} holds for $\mathcal{X}_1^{(p_1+1)}$ or
$\mathcal{X}_2^{(p_2)}$ then conclusion~(2) holds for $X^{(p)}$.  
Thus we may
assume that conclusion~(1) holds for both 
$\mathcal{X}_1^{(p_1+1)}$ and
$\mathcal{X}_2^{(p_2)}$.  Then  we obtain minimal genus Heegaard splittings
$C_1^1 \cup_{\s^1} C_2^1$ and $C_1^2 \cup_{\s^2} C_2^2$
respectively and Hopf-Haken annuli $A_1 \subset
\mathcal{X}_1^{(p_1+1)}$ and $A_2 \subset \mathcal{X}_2^{(p_2)}$
for $\s^1$ and $\s^2$, with $A_1$ connecting $\del X$ to $T$ and
$A_2$ connecting $T$ to $T_1$, and $A_1 \cap T = A_2 \cap T$. Then
by Lemma~\ref{lem:genus after amalgamation},
Propositions~\ref{prop:amalgamation of minimal genus: connected
case}, and~\ref{pro:combined haken annuli}, we see that conclusion
(1) of Theorem~\ref{thm:hopf-annulus} holds.

This completes the proof of Theorem~\ref{thm:hopf-annulus}.
\end{proof}

\section{Morimoto's Conjecture for m-small knots}
\label{sec:morimoto}

\begin{proof}[Proof of Theorem~\ref{thm:morimoto's conjecture}]
For the \lq\lq if\rq\rq part of Theorem~\ref{thm:morimoto's conjecture},
recall that  $g(E(\#_{i \in I} K_i)) \le \s_{i \in I} g(E(K_i))$ and
$g(E(\#_{i\not\in I} K_i)) \le \s_{i \not\in I} g(E(K_i))$    hold in general.
If $E(\#_{i \in I} K_i)$ admits a primitive meridian then by
Proposition~\ref{pro:primitive meridian implies degeneration} we have that
$g(X) < g(E(\#_{i \in I} K_i)) + g(E(\#_{i \in I} K_i))$.  Combining these
inequalities shows that $g(X) < \s_{i=1}^n g(E(K_i))$ as required.

Assume for a contradiction that the \lq\lq only if\rq\rq part of
Theorem~\ref{thm:morimoto's conjecture} is false and let $\{ K_i | i=1,\dots,n
\}$ be a collection of m-small knots that gives a counterexample, which
minimizes $n$ among all such counterexamples. (Note that we may assume the
uniqueness of prime decomposition of $\#_{i=1}^n K_i$ by Claim 1 of the proof
of  Theorem~\ref{thm:sft}.)  Let $X_i = E(K_i)$ and $X=E(\#_{i=1}^n K_i)$;
this and all other notation is kept as in the previous sections. Let $\s
\subset X$ be a minimal genus Heegaard surface. By Corollary~\ref{cor:sft},
$\s$ weakly reduces to a swallow follow torus, say $T$, giving the
decomposition:

$$X = X_I^{(1)} \cup_T X_J,$$

\noindent
where $I \cup J = \{1,\dots,n\}$, $I \neq \emptyset$, $J \neq \emptyset$, $I 
\cap J = \emptyset$, $X_I =E(\#_{i \in I} K_i)$, and
$X_J =E(\#_{i \in J} K_i)$.

By Corollary~\ref{cor:when can the genus go down}, we know that
either $g(X_I) = g(X_I^{(1)})$ or $g(X_I) = g(X_I^{(1)}) - 1$.  If
$g(X_I) = g(X_I^{(1)})$, then by Corollary~\ref{cor:primitive
meridian (hopf section)}, $X_I$ admits a primitive meridian and we
have a contradiction to our assumption that $\{K_i\}$ forms a
counterexample.  We may therefore assume that $g(X_I^{(1)}) - 1 =
g(X_I)$. By Corollary~\ref{cor:when can the genus go down}, $g(X)
= g(X_I) + g(X_J)$. Therefore by assumption we have that $g(X_I) +
g(X_J) = g(X) < \s_{i=1}^n g(X_i) = \s_{i \in I} g(X_i) + \s_{i
\in J} g(X_i)$. Therefore, either $g(X_I) <  \s_{i \in I} g(X_i)$
or $g(X_J) < \s_{i \in J} g(X_i)$ must hold (say the former); by the 
assumption of minimality on $n$, $X_I$ cannot give a
counterexample; therefore for some $\emptyset\neq I' \subset I$,
the exterior of $\#_{i \in I'} K_i$ admits a primitive meridian,
contradicting the assumption that the collection $\{ K_i |
i=1,\dots , n \}$ provides a counterexample.

This completes the proof of Theorem~\ref{thm:morimoto's
conjecture}.
\end{proof}

\nocite{rieck-sedgwick2}\nocite{rieck-sedgwick1}
% ----------------------------------------------------------------

\end{document}